\documentclass[A4j,11pt]{article}
\usepackage{amssymb,amsmath,amsthm,amscd,amsfonts}
\usepackage{amsmath,amssymb,amsthm,amscd,ascmac,mathrsfs,bm,type1cm,multirow,array,enumerate,mathrsfs,tabularx}
\usepackage{graphicx}
\usepackage{graphics}

\newcommand{\vv}{{\bm{v}}} 
\newcommand{\vw}{{\bm{w}}} 

\begin{document}

\title{{\bf Equifocal submanifolds\\
with non-flat section and\\
topological Tits buildings
}}
\author{{\bf Naoyuki Koike}}
\date{}
\maketitle

\begin{abstract}
From the Lytchak's result for polar foliations on an irreducible simply connected symmetric space $G/K$ of compact type and rank greater than one, 
we can derive that there exists no equifocal submanifold with non-flat section whose codimension is greater than two in the symmetric space $G/K$.  
In the first-half part of this paper, we give a new proof of this non-existence theorem.  The recipi of our new proof is as follows.  
Suppose that there exists an equifocal submanifold $M$ with non-flat section whose codimension is greater than two in an irreducible symmertric space $G/K$ 
of compact type and rank greater than one.  
We introduce the notion of a slice topology of $G/K$ associated to $M$.  We consider the universal covering $\pi:\widehat{G/K}\to G/K$ of the slice topological 
space $G/K$ and give $\widehat{G/K}$ the manifold structure and the Riemannian metric such that $\pi$ is a Riemannian submersion onto the symmetric space 
$G/K$.  First we show that a simplicial decomposition of the Riemannian manifold $\widehat{G/K}$ gives an irreducible topological Tits building of 
spherical type and rank greater than two.  By applying Burns-Spatzier's theorem to this topological Tits building, 
we show that the Riemannian manifold $\widehat{G/K}$ is homothetic to the unit sphere.  
Furthermore, from this fact, we show that 
$G/K$ is isometric to a sphere, a complex projective space or a quaternionic projective space.  
This contradicts that $G/K$ is of rank greater than one.  This is the recipi of our proof.  
In the second-half part, we estimate the codimension of $M$ from above by using the multiplicities of the roots of the root system of $G/K$.  
As its result, we can show that there exists no equifocal submanifold with non-flat section 
in some irreducible simply connected symmetric spaces of compact type.  
\end{abstract}

\vspace{0.5truecm}

\section{Introduction} 
In 1995, C. L. Terng and G. Thorbergsson (\cite{TT}) introduced the notion of an equifocal submanifold in a symmetric space $G/K$ 
as a compact immersed submanifold $M$ in $G/K$ satisfying the following conditions:

\vspace{0.2truecm}

(E-i) The normal holonomy group of $M$ is trivial;

(E-ii) For each $x\in M$, the normal umbrella $\Sigma_x:=\exp^{\perp}(T_x^{\perp}M)$ is a totally geodesic in $G/K$, where $\exp^{\perp}$ is 
the normal exponential map of $M$;

(E-iii) If $\widetilde{\vv}$ is a parallel normal vector field on $M$ such that $\exp^{\perp}(\widetilde{\vv}_{x_0})$ is 
a focal point of multiplicity $k$ for some $x_0\in M$, then $\exp^{\perp}(\widetilde{\vv}_x)$ also is a focal point of multiplicity $k$ 
for all $x\in M$;

(E-iv) For each $x\in M$, the induced metric on $\Sigma_x$ is flat.  

\vspace{0.2truecm}

\noindent
The totally geodesic submanifold $\Sigma_x$ in (E-ii) is called a {\it section of} $M$ {\it through} $x$.  
The condition (E-iii) is equivalent to the following condition:

\vspace{0.2truecm}

(E-iii') For each parallel unit normal vector field $\widetilde{\vv}$ of $M$, the end-point map $\eta_{\widetilde{\vv}}:M\to G/K$ is 
constant rank, where $\eta_{\widetilde{\vv}}$ is defined by $\eta_{\widetilde{\vv}}(x):=\exp^{\perp}(\widetilde{\vv}_x)$ ($x\in M$).  

\vspace{0.2truecm}

In 2004, M. M. Alexandrino (\cite{A1}) defined the notion of an equifocal submanifold in a general complete Riemannian manifold 
as a (not necessarily compact) immersed submanifold satisfying the above conditions (E-i), (E-ii) and (E-iii).  Furthermore, 
if the above condition (E-iv) also holds, he called the submanifold an {\it equifocal submanifold with flat section}.  
In \cite{A1}, \cite{A2}, \cite{AG} and \cite{AT}, the terminology ``equifocal submanifold'' was used in the sense of \cite{A1}.  
In this paper, we also shall use the terminology ``equifocal submanifold'' in the sense of \cite{A1}.  
In the sequel, we assume the following.  

\vspace{0.5truecm}

\noindent
{\bf Assumption.}\ We assume that all equifocal submanifolds are compact.  

\vspace{0.5truecm}

\noindent
{\it Remark 1.1.}\ (i)\ Let $M$ be an immersed submanifold in a symmetric space $G/K$.  If, $g_{\ast}^{-1}(T^{\perp}_{gK}M)$ is 
a Lie triple system of $\mathfrak p:=T_{eK}(G/K)(\subset\mathfrak g)$ for any $gK\in M$, then $M$ was said to 
{\it have the Lie triple systematic normal bundle} in \cite{Ko3}, where $\mathfrak g$ is the Lie alebra of $G$.  
The condition (E-ii) holds if and only if $M$ has Lie triple systematic normal bundle.  

(ii)\ The section $\Sigma_x$ meets the parallel submanifold $\eta_{\widetilde{\vv}}(M)$ of $M$ orthogonally, 
where $\widetilde{\vv}$ is any non-focal parallel normal vector field of $M$ (see Proposition 2.2 of \cite{HLO}).  

\vspace{0.5truecm}

For equifocal submanifolds with non-flat section, some open problems remain, for example the following.

\vspace{0.5truecm}

\noindent
{\bf Open Problem.} {\sl Does there exist no equifocal submanifold with non-flat section in an irreducible simply connected 
symmetric space of compact type and rank greater than one?}

\vspace{0.5truecm}

%

To solve this open problem, it is important to analyze the structure of the sections of the equifocal submanifold with non-flat section.  
In 2008, we \cite{Ko4} proved the following fact for the sections of an equifocal submanifold with non-flat section.  

\vspace{0.5truecm}

\noindent
{\bf Fact 1.} {\it The sections of an equifocal submanifold with non-flat section in an irreducible symmetric space of compact type are isometric to 
a sphere or a real projective space of constant curvature.}

\vspace{0.5truecm}

In 2014, A. Lytchak (\cite{L}) proved that polar foliations of codimension greater than two on an irreducible symmetric space of compact type and 
rank greater than one are hyperpolar, where we note that he treated also the case where the symmetric space is reducible.  
Since parallel submanifolds and focal submanifolds of an equifocal submanifold give a polar foliation, we can derive 
the following non-existence theorem of equifocal submanifolds with non-flat section from Lytchak's result.  

\vspace{0.5truecm}

\noindent
{\bf Theorem A.} {\it There eixsts no equifocal submanifold with non-flat section whose codimension is greater than two in an irreducible simply connected 
symmetric space $G/K$ of compact type and rank greater than one.}

\vspace{0.5truecm}

The recipi of his proof is as follows.  
Let $\mathfrak F$ be a polar foliation of codimension greater than two on an irreducible symmetric space $G/K$ of compact type and rank greater than one.  
He introduced a new metric called a {\it horizontal metric} on $G/K$ associated to $\mathfrak F$, where we note that the topology on the leaves of 
$\mathfrak F$ induced from the new metric are discrete.  Denote by $d^{\rm hor}$ this metric.  
Suppose that $\mathfrak F$ is not hyperpolar.  Let $\pi:\widehat{G/K}\to G/K$ be the universal covering of $(G/K,d^{\rm hor})$.  
He showed that a simplicial decomposition of $\widehat{G/K}$ gives an irreducible topological Tits building of spehrical type and rank greater than two 
by using the result for metric characterizations of spherical buildings in \cite{CL}.  
By applying Burns-Spatzier's theorem (see the next paragraph) to this topological Tits building, he showed that $\widehat{G/K}$ is homeomorphic to a sphere.  
Furthermore, from this fact, he showed that the fundamental group of $(G/K,d^{\rm hor})$ is a trivial group, the circle group $U(1)$ or the special unitary 
group $SU(2)$ and hence $G/K$ is homeomorphic to a sphere, a complex projective space or a quaternionic projective space.  This contradicts that $G/K$ is 
of rank greater than one.  Therefore, $\mathfrak F$ is hyperpolar.  This is the recipi of his proof.  

In 1987, K. Burns and R. Spatzier (\cite{BS}) classified irreducible topological Tits buildings of spherical type and rank greater than two 
such that the topology of the set of all vertices is metric topology.  In more detail, they showed that each of such irreducible topological 
Tits buildings of spherical type is isomorphic to the classical topological Tits building $(\mathcal B_{G'},\mathcal O_{G'})$ 
of spherical type associated to a simple non-compact (real) Lie group $G'$ without the center, where we note that $(\mathcal B_{G'},\mathcal O_{G'})$ is 
defined by using parabolic subgroups of $G'$.  See \cite{BS} about the detailed definition of $(\mathcal B_{G'},\mathcal O_{G'})$.  

In 1991, G. Thorbergsson (\cite{Th}) constructed the topological Tits building $(\mathcal B_M,\mathcal O_M)$ of spherical type associated to 
a full irreducible isoparametric submanifold $M$ of codimension greater than two in a Euclidean space.  As above, according to the classification theorem of 
Burns-Spatzier, the building $(\mathcal B_M,\mathcal O_M)$ is isomorphic to the classical topological Tits building $(\mathcal B_{G'},\mathcal O_{G'})$ 
associated to a simple non-compact (real) Lie group $G'$ without the center.  He constructed a homeomorphism of the ambient Euclidean space 
onto the tangent space of the symmetric space $G'/K$ of non-compact type, where $K'$ is the maximal compact Lie subgroup of $G'$.  
By showing that this homeomorphism is an isometry and the image of the isoparametric submanifold by this homeomorphism is a principal orbit of 
the isotropy representation of $G'/K$, he proved that the isoparametric submanifold is homogeneous.  

We give a new proof of Theorem A by using Burns-Spatzier's theorem and the discussion by Thorbergsson (\cite{Th}).  
The recipi of our new proof is as follows.  
Suppose that there exists an equifocal submanifold $M$ with non-flat section whose codimension is greater than two in an irreducible symmertric space $G/K$ 
of compact type and rank greater than one.  
We introduce a new topology of $G/K$ associated to $M$.  We call this new topology {\it slice topology} and denoted it by 
$\mathcal O_{\rm sl}$.  We consider the universal covering $\pi:\widehat{G/K}\to G/K$ of the topological space $(G/K,\mathcal O_{\rm sl})$ and give 
$\widehat{G/K}$ the manifold structure and the Riemannian metric such that $\pi$ is a Riemannian submersion onto the symmetric space $G/K$.  
First we show that a simplicial decomposition of the Riemannian manifold $\widehat{G/K}$ gives an irreducible topological Tits building of spherical type 
and rank greater than two by analyzing the structure of the simplicial decomposition.  
By applying Burns-Spatzier's theorem (see the next paragraph) to this topological Tits building, we show that the Riemannian manifold $\widehat{G/K}$ 
is homothetic to the unit sphere.  
Furthermore, from this fact, we show that the fundamental group of $(G/K,\mathcal O_{\rm sl})$ is a trivial group, the circle group $U(1)$ or 
the special unitary group $SU(2)$ and hence $G/K$ is isometric to a sphere, a complex projective space or a quaternionic projective space.  
This contradicts that $G/K$ is of rank greater than one.  Hence there does not exist the above equifocal submanifold.  This is the recipi of our proof.

In the second-half part of this paper, we estimate the codimension of $M$ from above by using the multiplicities of the roots of the root system of $G/K$.  
Let $\mathfrak a$ be a maximal abelian subspace of $\mathfrak p:=T_{eK}(G/K)\,(\subset\mathfrak g)$ and $\triangle_+$ the positive root system 
for $\mathfrak a$, where $e$ is the identity element of $G$ and $\mathfrak g$ is the Lie algebra of $G$.  
For each root $\alpha\in\triangle_+$, denote by $\mathfrak p_{\alpha}$ the root space for $\alpha$ and set 
$m_{\alpha}:={\rm dim}\,\mathfrak p_{\alpha}$.  

We prove the following fact for the estimate of the codimension of an equifocal submanifold with non-flat section.  

\vspace{0.5truecm}

\noindent
{\bf Theorem B.}\ {\sl Let $G/K$ be an irreducible simply connected symmetric space of compact type and $M$ an equifocal submanifold 
with non-flat section in $G/K$.  Then we have 
$${\rm codim}\,M\leq\left[\frac{1}{2}\,\mathop{\max}_{\alpha\in\triangle_+}(m_{\alpha}+m_{2\alpha})\right]+1,$$
where $m_{2\alpha}$ implies $0$ in the case of $2\alpha\notin\triangle_+$, and $[\cdot]$ is the floor function.  
}

\vspace{0.5truecm}

See Table 1 about the list of $\displaystyle{m_{G/K}:=\left[\frac{1}{2}\,\mathop{\max}_{\alpha\in\triangle_+}(m_{\alpha}+m_{2\alpha})\right]+1}$ 
for all irreducible simply connected symmetric spaces of compact type.  

\vspace{0.5truecm}

$$\begin{tabular}{|c|c|c|}
\hline
$G/K$ & $m_{G/K}$\\
\hline
$SU(n)/SO(n)\,\,\,\,(n\geq 3)$ & $1$\\
\hline
$SU(2n)/Sp(n)\,\,\,\,(n\geq3)$ & $3$\\
\hline
$SU(n)/S(U(i)\times U(n-i))\,\,\,\,(2\leq i<\frac{n}{2})$ & $n-2i+1$\\
\hline
$SU(2n)/S(U(n)\times U(n))\,\,\,\,(n\geq2)$ & $2$\\
\hline
$SO(2n+1)/(SO(i)\times SO(2n-i+1))\,\,\,\,(2\leq i\leq n)$ & $n-i+1$\\
\hline
$SO(2n)/(SO(i)\times SO(2n-i))\,\,\,\,(2\leq i\leq n-2)$ & $n-i+1$\\
\hline
$SO(2n)/(SO(n-1)\times SO(n+1))$ & $2$\\
\hline
$SO(2n)/(SO(n)\times SO(n))\,\,\,\,(n\geq2)$ & $1$\\
\hline
$SO(4n)/U((2n)\,\,\,\,(n\geq2)$ & $3$\\
\hline
$SO(4n+2)/U(2n+1)\,\,\,\,(n\geq2)$ & $3$\\
\hline
$Sp(n)/U(n)\,\,\,\,(n\geq2)$ & $1$\\
\hline
$Sp(n)/(Sp(i)\times Sp(n-i))\,\,\,\,(2\leq i<\frac{n}{2})$ & $2(n-2i+1)$\\
\hline
$Sp(2n)/(Sp(n)\times Sp(n))\,\,\,\,(n\geq2)$ & $3$\\
\hline
$E_6/Sp(4)$ & $1$\\
\hline
$E_6/SU(6)\cdot SU(2)$ & $2$\\
\hline
$E_6/Spin(10)\cdot U(1)$ & $5$\\
\hline
$E_6/F_4$ & $5$\\
\hline
$E_7/SU(8)$ & $1$\\
\hline
$E_7/SO(12)\cdot SU(2)$ & $3$\\
\hline
$E_7/E_6\cdot S^1$ & $5$\\
\hline
$E_8/SO(16)\cdot SU(2)$ & $1$\\
\hline
$E_8/E_7\cdot SU(2)$ & $5$\\
\hline
$F_4/Sp(3)\cdot SU(2)$ & $1$\\
\hline
$G_2/SO(4)$ & $1$\\
\hline
$(SU(n)\times SU(n))/\triangle SU(n)\,\,\,\,(n\geq3)$ & $2$\\
\hline
$(Sp(n)\times Sp(n))/\triangle Sp(n)\,\,\,\,(n\geq2)$ & $2$\\
\hline
$(E_6\times E_6)/\triangle E_6$ & $2$\\
\hline
$(E_7\times E_7)/\triangle E_7$ & $2$\\
\hline
$(E_8\times E_8)/\triangle E_8$ & $2$\\
\hline
$(F_4\times F_4)/\triangle F_4$ & $2$\\
\hline
$(G_2\times G_2)/\triangle G_2$ & $2$\\
\hline
\end{tabular}$$

\vspace{0.3truecm}

\centerline{{\bf Table 1: The list of $\displaystyle{m_{G/K}:=\left[\frac{1}{2}\mathop{\max}_{\alpha\in\triangle_+}(m_{\alpha}+m_{2\alpha})\right]+1}$}}

\newpage

From Table 1, we obtain the following fact.  

\vspace{0.5truecm}

\noindent
{\bf Theorem C.} 
{\sl There exists no equifocal submanifold with non-flat section in the following irreducible simply connected symmetric spaces of compact type:
{\small
$$\begin{array}{c}
SU(n)/SO(n)\,(n\geq 3),\,\,SO(2n)(SO(n)\times SO(n))\,(n\geq 2),\,\,
Sp(n)/U(n)\,(n\geq 2),\\
E_6/Sp(4),\,\,,E_7/SU(8),\,\,E_8/SO(16)\cdot SU(2),\,\,F_4/Sp(3)\cdot SU(2),\,\,
G_2/SO(4).
\end{array}$$
}
}

\section{Basic notions} 
In this section, we recall the basic notions.  First we recall some notions associated to an equifocal submanifold.  
Let $M$ be an $n$-dimensional equifocal submanifold with non-flat section in a symmetric space $G/K$ of compact type.  
Set $r:={\rm codim}\,M\,(\geq 2)$.  Denote by $\Sigma_x$ the section of $M$ through $x$.  
Let $\widetilde{\vv}$ be a parallel normal vector field of $M$ and $\eta_{\widetilde v}$ the end-point map for 
$\widetilde{\vv}$, which is the map of $M$ into $G/K$ defined by $\eta_{\widetilde{\vv}}(x):=\exp^{\perp}(\widetilde{\vv}_x)\,\,\,(x\in M)$, 
where $\exp^{\perp}$ is the normal exponential map of $M$.  Since $M$ is equifocal, $\eta_{\widetilde{\vv}}$ is of constant rank.  
If $\eta_{\widetilde{\vv}}$ is of constant rank smaller than $n$, then $\widetilde{\vv}$ is called a 
{\it focal normal vector field} of $M$ and $F_{\widetilde{\vv}}:=\eta_{\widetilde{\vv}}(M)$ is called the {\it focal submanifold} for 
$\widetilde{\vv}$.  Then $\eta_{\widetilde{\vv}}$ is a submersion of $M$ onto $F_{\widetilde{\vv}}$.  Each fibre of 
$\eta_{\widetilde{\vv}}$ ara called a {\it focal leaf} for $\widetilde{\vv}$.  Denote by $L^{\widetilde{\vv}}_x$ the focal leaf for 
$\widetilde{\vv}$ through $x\in M$.  
Define a distribution $\mathcal D_{\widetilde{\vv}}$ on $M$ by $(\mathcal D_{\widetilde{\vv}})_x:={\rm Ker}\,(d\eta_{\widetilde{\vv}})_x$ 
($x\in M$).  This distribution $\mathcal D_{\widetilde{\vv}}$ is called the {\it focal distribution} for $\widetilde{\vv}$.  
Set 
$${\mathcal TF}_x(M):=\{\widetilde{\vv}_x\,\vert\,\widetilde{\vv}\,:\,{\rm focal}\,\,{\rm normal}\,\,{\rm vector}\,\,{\rm field}
\,\,{\rm of}\,\,M\},$$
which is called the {\it tangential focal set} of $M$ at $x$.  Also, set ${\mathcal F}_x(M):=\exp^{\perp}({\mathcal TF}_x(M))$, 
which is called the {\it focal set} of $M$ at $x$.  
According to Fact 1 stated in Introduction, the section $\Sigma_x$ of $M$ through $x$ is isometric to an $r$-dimensional sphere 
(or an $r$-dimensional real projective space) of constant curvature.  
On the other hand, according to Lemma 2.15 and Proposition 2.16 of \cite{E}, the focal set $\mathcal F_x(M)$ consists of finitey many complete 
totally geodesic hypersurfaces in $\Sigma_x$ and they give a simplicial decomposition of $\Sigma_x$.  
We call the complete totally geodesic hypersurfaces in $\Sigma_x$ {\it focal walls of} $M$ {\it at} $x$.  
Denote by $\mathcal{FW}_x(M)$ the set of all the focal walls of $M$ at $x$.  
See Figure 1 about the graphical image of focal submanifolds, focal leaves and focal walls.  

\vspace{0.5truecm}

\centerline{
\unitlength 0.1in
\begin{picture}( 52.9400, 19.4600)( -1.5000,-25.5600)
%
\special{pn 8}%
\special{ar 4160 2294 754 864  3.2028299 4.1791740}%
%
\special{pn 8}%
\special{ar 5416 2380 754 864  3.2027853 4.1783186}%
%
\special{pn 8}%
\special{ar 3394 2700 2136 470  4.7186972 5.3528522}%
%
\special{pn 8}%
\special{ar 3756 2026 2136 470  4.7187041 5.3540783}%
%
\special{pn 8}%
\special{ar 1684 2162 756 864  3.2035570 4.1803242}%
%
\special{pn 8}%
\special{ar 2940 2248 756 864  3.2021026 4.1798978}%
%
\special{pn 8}%
\special{ar 920 2568 2134 470  4.7187076 5.3540127}%
%
\special{pn 8}%
\special{ar 1282 1896 2134 468  4.7186868 5.3523394}%
%
\special{pn 13}%
\special{ar 3424 2478 2136 470  4.7187128 5.3545943}%
%
\special{pn 13}%
\special{ar 4612 2620 1198 1090  3.5081312 4.4626371}%
%
\special{pn 13}%
\special{ar 3668 2610 944 1206  5.2414657 6.0333175}%
%
\special{pn 20}%
\special{sh 1}%
\special{ar 1662 1776 10 10 0  6.28318530717959E+0000}%
\special{sh 1}%
\special{ar 1662 1776 10 10 0  6.28318530717959E+0000}%
%
\special{pn 8}%
\special{ar 2574 2224 526 1436  3.2015392 4.0680109}%
%
\special{pn 8}%
\special{ar 2574 2224 516 1426  2.9747157 3.1415927}%
%
\special{pn 8}%
\special{ar 1770 2162 524 1436  3.2015392 4.0691804}%
%
\special{pn 8}%
\special{ar 1752 2154 518 1426  2.9745992 3.1415927}%
%
\special{pn 13}%
\special{ar 2164 2180 526 1434  3.2014840 4.0688879}%
%
\special{pn 13}%
\special{ar 2156 2172 516 1426  2.9747157 3.1415927}%
%
\special{pn 13}%
\special{ar 1002 2114 2122 350  4.7187460 5.3534160}%
%
\special{pn 8}%
\special{pa 1508 754}%
\special{pa 1832 1064}%
\special{da 0.070}%
\special{sh 1}%
\special{pa 1832 1064}%
\special{pa 1798 1004}%
\special{pa 1794 1028}%
\special{pa 1770 1032}%
\special{pa 1832 1064}%
\special{fp}%
%
\special{pn 8}%
\special{pa 1250 890}%
\special{pa 1528 1112}%
\special{da 0.070}%
\special{sh 1}%
\special{pa 1528 1112}%
\special{pa 1488 1054}%
\special{pa 1486 1078}%
\special{pa 1464 1086}%
\special{pa 1528 1112}%
\special{fp}%
%
\special{pn 8}%
\special{pa 976 1498}%
\special{pa 1106 1764}%
\special{da 0.070}%
\special{sh 1}%
\special{pa 1106 1764}%
\special{pa 1094 1694}%
\special{pa 1082 1716}%
\special{pa 1060 1712}%
\special{pa 1106 1764}%
\special{fp}%
%
\special{pn 20}%
\special{sh 1}%
\special{ar 2076 1804 10 10 0  6.28318530717959E+0000}%
\special{sh 1}%
\special{ar 2076 1804 10 10 0  6.28318530717959E+0000}%
%
\special{pn 8}%
\special{pa 1092 1310}%
\special{pa 1380 1656}%
\special{da 0.070}%
\special{sh 1}%
\special{pa 1380 1656}%
\special{pa 1352 1592}%
\special{pa 1346 1614}%
\special{pa 1322 1618}%
\special{pa 1380 1656}%
\special{fp}%
\put(8.7000,-14.8000){\makebox(0,0)[lb]{$\Sigma_x$}}%
\put(9.8300,-12.8100){\makebox(0,0)[lb]{$L^{\widetilde v}_x$}}%
\put(12.1900,-9.1900){\makebox(0,0)[rb]{$M$}}%
\put(14.7000,-7.8000){\makebox(0,0)[rb]{$F_{\widetilde v}$}}%
\put(20.9500,-18.4200){\makebox(0,0)[lt]{$x$}}%
\put(34.0300,-24.0400){\makebox(0,0)[lt]{{\small focal walls of $M$ at $x$}}}%
%
\special{pn 8}%
\special{ar 2950 2214 898 712  3.7447902 3.8193243}%
\special{ar 2950 2214 898 712  3.8640448 3.9385790}%
\special{ar 2950 2214 898 712  3.9832995 4.0578336}%
\special{ar 2950 2214 898 712  4.1025541 4.1770883}%
\special{ar 2950 2214 898 712  4.2218088 4.2963430}%
\special{ar 2950 2214 898 712  4.3410635 4.4155976}%
\special{ar 2950 2214 898 712  4.4603181 4.5348523}%
\special{ar 2950 2214 898 712  4.5795728 4.6541069}%
\special{ar 2950 2214 898 712  4.6988274 4.7733616}%
\special{ar 2950 2214 898 712  4.8180821 4.8926163}%
\special{ar 2950 2214 898 712  4.9373367 5.0118709}%
\special{ar 2950 2214 898 712  5.0565914 5.1311256}%
\special{ar 2950 2214 898 712  5.1758461 5.2298878}%
%
\special{pn 8}%
\special{pa 3394 1604}%
\special{pa 3458 1634}%
\special{fp}%
\special{sh 1}%
\special{pa 3458 1634}%
\special{pa 3406 1588}%
\special{pa 3410 1612}%
\special{pa 3390 1624}%
\special{pa 3458 1634}%
\special{fp}%
\put(28.2200,-14.4200){\makebox(0,0)[lb]{in fact}}%
\put(22.0300,-9.0000){\makebox(0,0)[lb]{$\exp^{\perp}(T^{\perp}_{\eta_{\widetilde v}(x)}F_{\widetilde v})$}}%
%
\special{pn 8}%
\special{pa 2528 938}%
\special{pa 2242 1642}%
\special{da 0.070}%
\special{sh 1}%
\special{pa 2242 1642}%
\special{pa 2286 1588}%
\special{pa 2262 1594}%
\special{pa 2248 1574}%
\special{pa 2242 1642}%
\special{fp}%
%
\special{pn 20}%
\special{sh 1}%
\special{ar 4230 1918 10 10 0  6.28318530717959E+0000}%
\special{sh 1}%
\special{ar 4230 1918 10 10 0  6.28318530717959E+0000}%
%
\special{pn 20}%
\special{sh 1}%
\special{ar 3866 1786 10 10 0  6.28318530717959E+0000}%
\special{sh 1}%
\special{ar 3866 1786 10 10 0  6.28318530717959E+0000}%
\put(42.6800,-18.8000){\makebox(0,0)[lt]{$x$}}%
\put(37.6700,-14.7100){\makebox(0,0)[lb]{{\small $\eta_{\widetilde v}(x)$}}}%
\put(41.0100,-18.1400){\makebox(0,0)[lb]{$\widetilde v_x$}}%
%
\special{pn 8}%
\special{ar 3580 2594 1112 838  4.9597558 5.3350598}%
%
\special{pn 8}%
\special{pa 3954 1490}%
\special{pa 3856 1766}%
\special{da 0.070}%
\special{sh 1}%
\special{pa 3856 1766}%
\special{pa 3896 1710}%
\special{pa 3874 1716}%
\special{pa 3858 1696}%
\special{pa 3856 1766}%
\special{fp}%
\put(37.3700,-25.5600){\makebox(0,0)[lt]{$\,$}}%
%
\special{pn 8}%
\special{pa 4220 1910}%
\special{pa 3934 1748}%
\special{fp}%
\special{sh 1}%
\special{pa 3934 1748}%
\special{pa 3982 1798}%
\special{pa 3980 1774}%
\special{pa 4002 1764}%
\special{pa 3934 1748}%
\special{fp}%
%
\special{pn 13}%
\special{pa 2066 1814}%
\special{pa 2054 1844}%
\special{pa 2036 1870}%
\special{pa 2008 1888}%
\special{pa 1980 1902}%
\special{pa 1950 1914}%
\special{pa 1920 1924}%
\special{pa 1890 1932}%
\special{pa 1858 1936}%
\special{pa 1826 1942}%
\special{pa 1794 1944}%
\special{pa 1762 1946}%
\special{pa 1730 1948}%
\special{pa 1698 1946}%
\special{pa 1666 1946}%
\special{pa 1634 1942}%
\special{pa 1602 1940}%
\special{pa 1570 1936}%
\special{pa 1540 1932}%
\special{pa 1508 1924}%
\special{pa 1478 1916}%
\special{pa 1446 1906}%
\special{pa 1416 1896}%
\special{pa 1388 1882}%
\special{pa 1358 1868}%
\special{pa 1332 1852}%
\special{pa 1308 1830}%
\special{pa 1288 1804}%
\special{pa 1278 1774}%
\special{pa 1278 1766}%
\special{sp}%
%
\special{pn 13}%
\special{ar 1672 1758 392 142  3.1415927 4.6396779}%
%
\special{pn 13}%
\special{pa 1732 1624}%
\special{pa 1764 1626}%
\special{pa 1796 1630}%
\special{pa 1826 1636}%
\special{pa 1858 1640}%
\special{pa 1890 1648}%
\special{pa 1920 1654}%
\special{pa 1950 1664}%
\special{pa 1980 1676}%
\special{pa 2010 1688}%
\special{pa 2038 1704}%
\special{pa 2062 1724}%
\special{pa 2080 1750}%
\special{pa 2084 1782}%
\special{pa 2084 1786}%
\special{sp}%
%
\special{pn 13}%
\special{pa 2164 1338}%
\special{pa 2154 1368}%
\special{pa 2136 1394}%
\special{pa 2110 1412}%
\special{pa 2082 1430}%
\special{pa 2054 1442}%
\special{pa 2024 1454}%
\special{pa 1994 1462}%
\special{pa 1962 1468}%
\special{pa 1930 1474}%
\special{pa 1898 1478}%
\special{pa 1866 1480}%
\special{pa 1834 1480}%
\special{pa 1802 1480}%
\special{pa 1770 1480}%
\special{pa 1738 1478}%
\special{pa 1706 1476}%
\special{pa 1676 1470}%
\special{pa 1644 1464}%
\special{pa 1612 1460}%
\special{pa 1582 1450}%
\special{pa 1550 1442}%
\special{pa 1520 1432}%
\special{pa 1492 1418}%
\special{pa 1462 1404}%
\special{pa 1436 1388}%
\special{pa 1410 1368}%
\special{pa 1388 1346}%
\special{pa 1374 1316}%
\special{pa 1368 1286}%
\special{pa 1368 1286}%
\special{sp}%
%
\special{pn 13}%
\special{pa 1368 1310}%
\special{pa 1376 1280}%
\special{pa 1396 1254}%
\special{pa 1422 1236}%
\special{pa 1450 1220}%
\special{pa 1478 1206}%
\special{pa 1508 1196}%
\special{pa 1540 1188}%
\special{pa 1572 1182}%
\special{pa 1602 1178}%
\special{pa 1634 1172}%
\special{pa 1666 1170}%
\special{pa 1698 1168}%
\special{pa 1730 1168}%
\special{pa 1760 1166}%
\special{sp}%
%
\special{pn 13}%
\special{ar 1790 1320 374 144  4.8644791 6.2831853}%
\special{ar 1790 1320 374 144  0.0000000 0.0862056}%
%
\special{pn 13}%
\special{ar 3708 1996 2118 396  4.7184579 5.3522438}%
%
\special{pn 13}%
\special{ar 5450 2692 1252 946  3.5916883 4.1842243}%
%
\special{pn 13}%
\special{ar 3050 2624 748 952  5.3686399 5.8709527}%
\put(51.4400,-20.9000){\makebox(0,0)[lb]{$\Sigma_x$}}%
\put(27.0400,-12.0400){\makebox(0,0)[lb]{{\small $\eta_{\widetilde v}(x)$}}}%
%
\special{pn 8}%
\special{pa 2646 1214}%
\special{pa 1672 1766}%
\special{da 0.070}%
\special{sh 1}%
\special{pa 1672 1766}%
\special{pa 1740 1750}%
\special{pa 1718 1740}%
\special{pa 1720 1716}%
\special{pa 1672 1766}%
\special{fp}%
\end{picture}%
\hspace{3truecm}}

\vspace{0.5truecm}

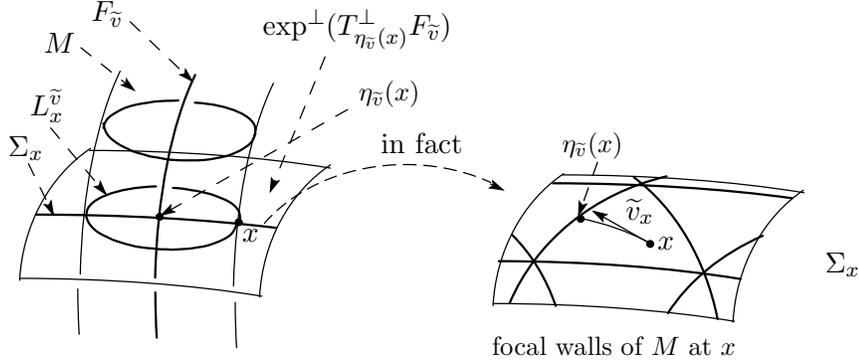
\centerline{{\bf Figure 1$\,:\,$ Focal submanifold, focal leaf and focal wall}}

\vspace{0.5truecm}

Next we recall the notion of a topological Tits building.  
Let $\Delta=({\mathcal V},{\mathcal S})$ be an $r$-dimensional simplicial complex, where 
${\mathcal V}$ denotes the set of all vertices and ${\mathcal S}$ denotes the set of all simplices.  
Each $r$-simplex of $\Delta$ is called a {\it chamber} of $\Delta$.  
Let ${\mathcal A}:=\{{\mathcal A}_{\lambda}\}_{\lambda\in\Lambda}$ be a family of subcomplexes of $\Delta$.  
The pair $\mathcal B:=(\Delta,\mathcal A)$ is called a {\it Tits building} if the following conditions hold:

\vspace{0.2truecm}

\noindent
(B1) Each $(r-1)$-dimensional simplex of $\Delta$ is contained in at least three chambers.  

\noindent
(B2) Each $(r-1)$-dimensional simplex in a subcomplex ${\mathcal A}_{\lambda}$ are contained in exactly two chambers 
of ${\mathcal A}_{\lambda}$.  

\noindent
(B3) Any two simplices of $\Delta$ are contained in some ${\mathcal A}_{\lambda}$.  

\noindent
(B4) If two subcomplexes ${\mathcal A}_{\lambda_1}$ and ${\mathcal A}_{\lambda_2}$ share a chamber, then there is 
an isomorphism of ${\mathcal A}_{\lambda_1}$ onto ${\mathcal A}_{\lambda_2}$ fixing 
${\mathcal A}_{\lambda_1}\cap{\mathcal A}_{\lambda_2}$ pointwisely.  

\vspace{0.2truecm}

\noindent 
Each subcomplex belonging to ${\mathcal A}$ is called an {\it apartment} of $\mathcal B$.  
In this paper, we assume that all Tits building furthermore satisfies the following condition: 

\vspace{0.2truecm}

\noindent
(B5) Each apartment ${\mathcal A}_{\lambda}$ is a coxeter complex.  

\vspace{0.2truecm}

\noindent
If each ${\mathcal A}_{\lambda}$ is finite (resp. infinite), then the building ${\mathcal B}$ is said to be {\it spherical type} 
(resp. {\it affine type}).  Let ${\mathcal O}$ be a Hausdorff topology of ${\mathcal V}$.  The pair $({\mathcal B},{\mathcal O})$ 
is called a {\it topological Tits building} if the following conditions hold:

\vspace{0.2truecm}

\noindent
(TB1) $\mathcal B$ is a Tits building.  

\noindent
(TB2) For $k\in\{1,\cdots,r\}$, $\widehat{\mathcal S}_k:=\{(x_1,\cdots,x_{k+1})\in{\mathcal V}^{k+1}\,\vert\,\,
\vert x_1\cdots x_{k+1}\vert\in{\mathcal S}_k\}$ is closed in the product topological space 
$({\mathcal V}^{k+1},{\mathcal O}^{k+1})$, where ${\mathcal S}_k$ denotes the set of all $k$-simplices of 
${\mathcal S}$ and $\vert x_1\cdots x_{k+1}\vert$ denotes the $k$-simplex with vertices $x_1,\cdots,x_{k+1}$.  

\vspace{0.2truecm}

\noindent
An $r$-dimensional simplical complex $\Delta$ is called a {\it chamber complex} if, for any two chambers $C,\,C'$ of $\Delta$, 
there exists a sequence $\{C_1=C,C_2,\cdots,C_k=C'\}$ such that $C_i\cap C_{i+1}$ ($i=1,\cdots,k-1$) are $(r-1)$-dimensional simplex.  

\section{A new proof of Theorem A (spherical section-case)} 
Let $G/K$ be an irreducible simply connected symmetric space of compact type and rank greater than one, and $M$ an equifocal submanifold with non-flat section 
in $G/K$.  We consider the case where the sections are spheres.  
Now we shall construct a topological Tits building of spherical type associated to $M$ defined on $G/K$.  
Denote by $r$ the codimension of $M$.  
Let $\iota$ be the natural embedding of the sphere $\Sigma_x$ into an $(r+1)$-dimensional Euclidean space $\mathbb R^{r+1}$.  
Since each element $S$ of $\mathcal{FW}_x(M)$ is a totally geodesic hypersphere in $\Sigma_x$, the cone over $\iota(S)$ is an $r$-dimensional 
vector subspace of $\mathbb R^{r+1}$.  Denote by $\Pi_S$ this $r$-dimensional vector subspace.  
Since $G/K$ is irreducible, $M$ is a full irreducible submanifold.  Hence the group generated by the reflections with respect to $\Pi_S$ 
is an irreducible finite coxeter group of rank $r$, that is, $\mathcal{FW}_x(M)$ gives a coxeter complex, which has a $(r+1)$-simplex as the chambers, 
defined on the sphere $\Sigma_x$.  Denote by $\mathcal A_x=(\mathcal V_x,\mathcal S_x)$ this coxeter complex defined on $\Sigma_x$.  
Set $\displaystyle{\mathcal S_M:=\mathop{\cup}_{x\in M}\mathcal S_x}$ and $\displaystyle{\mathcal V_M:=\mathop{\cup}_{x\in M}\mathcal V_x}$.  
It is clear that $\Delta_M:=(\mathcal V_M,\mathcal S_M)$ is a simplicial complex with $|\Delta_M|=G/K$.  
When $\mathcal A_{x_1}\cap\mathcal A_{x_2}$ contains a face of a chamber (i.e., $(r-1)$-dimensional simplex) of $\Delta_M$, 
we can construct four coxeter subcomplexes of $\Delta_M$ by patching subcomplexes of $\mathcal A_{x_1}$ and $\mathcal A_{x_2}$ (see Figure 2).  
Denote by $\mathcal A_{x_1x_2}$ one of the four coxeter subcomplexes of $\Delta_M$.  Furthermore, when 
$\mathcal A_{x_1x_2}\cap\mathcal A_{x_3}$ contains a face of a chamber (i.e., $(r-1)$-dimensional simplex) of $\Delta_M$, we can construct 
four coxeter subcomplexes by patching subcomplexes of $\mathcal A_{x_1x_2}$ and $\mathcal A_{x_3}$.  
In the sequel, we can construct infinitely many coxeter subcomplexes of $\Delta_M$ by repeating this process.  
Let $\mathcal A_M=\{\mathcal A_{\lambda}\}_{\lambda\in\Lambda}$ be the sum of $\{\mathcal A_x\}_{x\in M}$ and the set of all coxeter 
subcomplexes of $\Delta_M$ constructed thus.  Set $\mathcal B_M:=(\Delta_M,\mathcal A_M)$.  It is clear that $\mathcal V_M$ is equal to 
the sum of $(r+1)$-pieces of focal submanifolds of $M$ (see Figure 3).  Give $\mathcal V_M$ the topology induced from the topology of $G/K$.  
Denote by $\mathcal O_M$ this topology.  

\vspace{0.35truecm}

\centerline{
\unitlength 0.1in
\begin{picture}( 19.7900, 20.3700)( 26.7600,-27.9200)
\put(42.1300,-27.9200){\makebox(0,0)[lt]{$\,$}}%
%
\special{pn 8}%
\special{ar 3398 2222 2044 428  4.7190019 5.3542665}%
%
\special{pn 13}%
\special{ar 2804 3078 2500 910  4.7184836 5.3535704}%
%
\special{pn 13}%
\special{ar 4222 2766 1144 990  3.5064977 4.4621232}%
%
\special{pn 13}%
\special{ar 3258 2802 1016 1148  5.2400921 6.0327563}%
%
\special{pn 20}%
\special{sh 1}%
\special{ar 3800 2118 10 10 0  6.28318530717959E+0000}%
\special{sh 1}%
\special{ar 3800 2118 10 10 0  6.28318530717959E+0000}%
\put(38.5500,-20.4000){\makebox(0,0)[lt]{{\small $x_1$}}}%
\put(46.5500,-17.1100){\makebox(0,0)[lb]{$\mathcal A_{x_1}$}}%
%
\special{pn 20}%
\special{sh 1}%
\special{ar 3290 2196 10 10 0  6.28318530717959E+0000}%
\special{sh 1}%
\special{ar 3290 2196 10 10 0  6.28318530717959E+0000}%
%
\special{pn 13}%
\special{ar 2004 1996 1838 1428  5.4932798 6.2831853}%
%
\special{pn 8}%
\special{ar 3000 1936 960 1886  5.6168443 6.2776537}%
%
\special{pn 8}%
\special{ar 2150 2728 1008 2416  5.6160665 6.2779129}%
%
\special{pn 8}%
\special{ar 4994 3044 1938 1342  4.1308172 4.1470485}%
%
\special{pn 8}%
\special{ar 5370 2066 2832 1600  3.6901950 4.1018890}%
\put(41.1900,-9.5000){\makebox(0,0)[lb]{$\mathcal A_{x_2}$}}%
\put(46.2700,-12.8800){\makebox(0,0)[lb]{$\mathcal A_{x_1x_2}$}}%
%
\special{pn 8}%
\special{pa 4090 906}%
\special{pa 3696 1062}%
\special{da 0.070}%
\special{sh 1}%
\special{pa 3696 1062}%
\special{pa 3766 1056}%
\special{pa 3746 1042}%
\special{pa 3752 1020}%
\special{pa 3696 1062}%
\special{fp}%
%
\special{pn 8}%
\special{ar 4636 2048 302 408  3.3314092 3.5004233}%
\special{ar 4636 2048 302 408  3.6018318 3.7708458}%
\special{ar 4636 2048 302 408  3.8722543 4.0412684}%
\special{ar 4636 2048 302 408  4.1426768 4.3116909}%
\special{ar 4636 2048 302 408  4.4130994 4.5821134}%
\special{ar 4636 2048 302 408  4.6835219 4.7123890}%
%
\special{pn 8}%
\special{pa 4354 1918}%
\special{pa 4346 1962}%
\special{da 0.070}%
\special{sh 1}%
\special{pa 4346 1962}%
\special{pa 4378 1902}%
\special{pa 4356 1910}%
\special{pa 4340 1894}%
\special{pa 4346 1962}%
\special{fp}%
%
\special{pn 8}%
\special{pa 4582 1244}%
\special{pa 3810 1382}%
\special{da 0.070}%
\special{sh 1}%
\special{pa 3810 1382}%
\special{pa 3878 1390}%
\special{pa 3862 1374}%
\special{pa 3872 1352}%
\special{pa 3810 1382}%
\special{fp}%
%
\special{pn 8}%
\special{ar 4636 1980 518 736  3.1954763 3.2912465}%
\special{ar 4636 1980 518 736  3.3487085 3.4444787}%
\special{ar 4636 1980 518 736  3.5019408 3.5977109}%
\special{ar 4636 1980 518 736  3.6551730 3.7509432}%
\special{ar 4636 1980 518 736  3.8084053 3.9041754}%
\special{ar 4636 1980 518 736  3.9616375 4.0574077}%
\special{ar 4636 1980 518 736  4.1148698 4.2106399}%
\special{ar 4636 1980 518 736  4.2681020 4.3638722}%
\special{ar 4636 1980 518 736  4.4213342 4.5171044}%
\special{ar 4636 1980 518 736  4.5745665 4.6064651}%
%
\special{pn 8}%
\special{pa 4120 1936}%
\special{pa 4120 1970}%
\special{da 0.070}%
\special{sh 1}%
\special{pa 4120 1970}%
\special{pa 4140 1904}%
\special{pa 4120 1918}%
\special{pa 4100 1904}%
\special{pa 4120 1970}%
\special{fp}%
%
\special{pn 8}%
\special{ar 4928 3268 1892 1548  3.5026524 3.7032009}%
%
\special{pn 8}%
\special{ar 4844 3174 1770 1420  3.7104965 3.8593824}%
%
\special{pn 8}%
\special{ar 4844 3174 1750 1454  3.8900439 4.1340473}%
%
\special{pn 8}%
\special{ar 5042 2906 846 1212  3.3596983 4.1648267}%
%
\special{pn 8}%
\special{ar 3650 3088 1120 1412  4.2109516 4.3169087}%
%
\special{pn 8}%
\special{ar 3462 2818 846 1058  3.5410039 4.1591725}%
%
\special{pn 8}%
\special{ar 3058 2482 2964 694  4.7654626 4.8040519}%
%
\special{pn 8}%
\special{ar 2512 2922 2042 528  4.7934091 5.6923894}%
%
\special{pn 4}%
\special{pa 3176 2410}%
\special{pa 3150 2386}%
\special{dt 0.027}%
\special{pa 3246 2422}%
\special{pa 3142 2328}%
\special{dt 0.027}%
\special{pa 3316 2434}%
\special{pa 3146 2278}%
\special{dt 0.027}%
\special{pa 3386 2448}%
\special{pa 3146 2226}%
\special{dt 0.027}%
\special{pa 3456 2460}%
\special{pa 3136 2166}%
\special{dt 0.027}%
\special{pa 3526 2472}%
\special{pa 3122 2100}%
\special{dt 0.027}%
\special{pa 3596 2484}%
\special{pa 3122 2048}%
\special{dt 0.027}%
\special{pa 3666 2496}%
\special{pa 3116 1992}%
\special{dt 0.027}%
\special{pa 3736 2510}%
\special{pa 3110 1934}%
\special{dt 0.027}%
\special{pa 3806 2522}%
\special{pa 3104 1878}%
\special{dt 0.027}%
\special{pa 3876 2534}%
\special{pa 3090 1812}%
\special{dt 0.027}%
\special{pa 3946 2546}%
\special{pa 3074 1746}%
\special{dt 0.027}%
\special{pa 4016 2558}%
\special{pa 3062 1684}%
\special{dt 0.027}%
\special{pa 4086 2572}%
\special{pa 3058 1624}%
\special{dt 0.027}%
\special{pa 4156 2584}%
\special{pa 3040 1558}%
\special{dt 0.027}%
\special{pa 4226 2596}%
\special{pa 3024 1492}%
\special{dt 0.027}%
\special{pa 4232 2550}%
\special{pa 3010 1426}%
\special{dt 0.027}%
\special{pa 4242 2508}%
\special{pa 2988 1352}%
\special{dt 0.027}%
\special{pa 4252 2466}%
\special{pa 2964 1280}%
\special{dt 0.027}%
\special{pa 4264 2424}%
\special{pa 2968 1232}%
\special{dt 0.027}%
\special{pa 4278 2386}%
\special{pa 2998 1208}%
\special{dt 0.027}%
\special{pa 4294 2348}%
\special{pa 3028 1184}%
\special{dt 0.027}%
\special{pa 4306 2306}%
\special{pa 3058 1158}%
\special{dt 0.027}%
\special{pa 4320 2268}%
\special{pa 3088 1134}%
\special{dt 0.027}%
\special{pa 4336 2230}%
\special{pa 3120 1112}%
\special{dt 0.027}%
\special{pa 4352 2194}%
\special{pa 3950 1824}%
\special{dt 0.027}%
\special{pa 3950 1824}%
\special{pa 3152 1088}%
\special{dt 0.027}%
\special{pa 4376 2162}%
\special{pa 4006 1824}%
\special{dt 0.027}%
\special{pa 3950 1772}%
\special{pa 3182 1066}%
\special{dt 0.027}%
\special{pa 4394 2128}%
\special{pa 4072 1830}%
\special{dt 0.027}%
%
\special{pn 4}%
\special{pa 3948 1718}%
\special{pa 3214 1044}%
\special{dt 0.027}%
\special{pa 4412 2094}%
\special{pa 4132 1836}%
\special{dt 0.027}%
\special{pa 3942 1662}%
\special{pa 3246 1020}%
\special{dt 0.027}%
\special{pa 4432 2060}%
\special{pa 4194 1842}%
\special{dt 0.027}%
\special{pa 3938 1606}%
\special{pa 3278 996}%
\special{dt 0.027}%
\special{pa 4458 2032}%
\special{pa 4256 1844}%
\special{dt 0.027}%
\special{pa 3932 1548}%
\special{pa 3314 980}%
\special{dt 0.027}%
\special{pa 4490 2008}%
\special{pa 4316 1848}%
\special{dt 0.027}%
\special{pa 3926 1490}%
\special{pa 3350 962}%
\special{dt 0.027}%
\special{pa 4516 1982}%
\special{pa 4376 1854}%
\special{dt 0.027}%
\special{pa 3920 1432}%
\special{pa 3386 942}%
\special{dt 0.027}%
\special{pa 4540 1950}%
\special{pa 4440 1858}%
\special{dt 0.027}%
\special{pa 3914 1376}%
\special{pa 3420 920}%
\special{dt 0.027}%
\special{pa 4566 1924}%
\special{pa 4512 1874}%
\special{dt 0.027}%
\special{pa 3900 1312}%
\special{pa 3454 902}%
\special{dt 0.027}%
\special{pa 4596 1900}%
\special{pa 4578 1880}%
\special{dt 0.027}%
\special{pa 3886 1246}%
\special{pa 3490 882}%
\special{dt 0.027}%
\special{pa 3870 1178}%
\special{pa 3528 864}%
\special{dt 0.027}%
\special{pa 3852 1110}%
\special{pa 3562 842}%
\special{dt 0.027}%
\special{pa 3832 1042}%
\special{pa 3594 822}%
\special{dt 0.027}%
\special{pa 3814 970}%
\special{pa 3634 806}%
\special{dt 0.027}%
\special{pa 3792 900}%
\special{pa 3678 794}%
\special{dt 0.027}%
\special{pa 3770 826}%
\special{pa 3718 778}%
\special{dt 0.027}%
\put(36.0100,-16.3400){\makebox(0,0)[lt]{{\small $x_2$}}}%
%
\special{pn 20}%
\special{sh 1}%
\special{ar 3566 1686 10 10 0  6.28318530717959E+0000}%
\special{sh 1}%
\special{ar 3566 1686 10 10 0  6.28318530717959E+0000}%
%
\special{pn 13}%
\special{ar 4712 2516 1430 2380  3.1842148 3.9671366}%
%
\special{pn 13}%
\special{ar 4712 2516 1430 2362  3.1415927 3.1692536}%
\end{picture}%
\hspace{0.5truecm}}

\vspace{0.35truecm}

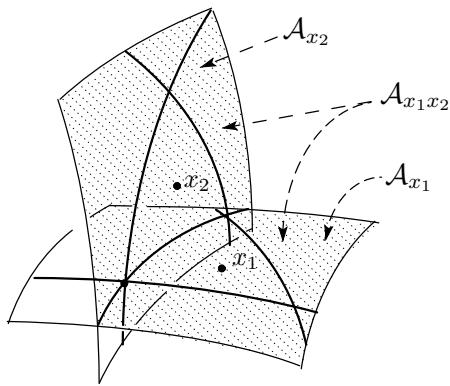
\centerline{{\bf Figure 2$\,:\,$ The coxeter complexes $\mathcal A_{x_1},\,\mathcal A_{x_2}$ and $\mathcal A_{x_1x_2}$}}

\vspace{0.5truecm}

\centerline{
\unitlength 0.1in
\begin{picture}( 45.5400, 24.9000)( 12.7400,-38.3900)
%
\special{pn 8}%
\special{ar 6170 3688 728 854  3.2026487 4.1785444}%
%
\special{pn 8}%
\special{ar 4214 4004 2066 464  4.7191877 5.3527507}%
%
\special{pn 8}%
\special{ar 4564 3338 2064 464  4.7186442 5.3551530}%
%
\special{pn 13}%
\special{ar 4284 3772 2036 468  4.7190133 5.3546676}%
%
\special{pn 13}%
\special{ar 5392 3924 1158 1076  3.5082189 4.4627680}%
%
\special{pn 13}%
\special{ar 4482 3916 912 1190  5.2414481 6.0331494}%
%
\special{pn 20}%
\special{sh 1}%
\special{ar 1976 2912 10 10 0  6.28318530717959E+0000}%
\special{sh 1}%
\special{ar 1976 2912 10 10 0  6.28318530717959E+0000}%
%
\special{pn 8}%
\special{ar 2856 3340 508 1418  3.2016349 4.0685855}%
%
\special{pn 8}%
\special{ar 2078 3280 508 1418  3.2016921 4.0685855}%
%
\special{pn 13}%
\special{ar 2458 3296 508 1416  3.2006279 4.0673773}%
%
\special{pn 13}%
\special{ar 1334 3232 2052 344  4.7186746 5.3524599}%
%
\special{pn 8}%
\special{pa 1574 2024}%
\special{pa 1844 2242}%
\special{da 0.070}%
\special{sh 1}%
\special{pa 1844 2242}%
\special{pa 1804 2186}%
\special{pa 1802 2208}%
\special{pa 1780 2216}%
\special{pa 1844 2242}%
\special{fp}%
%
\special{pn 20}%
\special{sh 1}%
\special{ar 2372 2926 10 10 0  6.28318530717959E+0000}%
\special{sh 1}%
\special{ar 2372 2926 10 10 0  6.28318530717959E+0000}%
\put(26.2900,-28.9800){\makebox(0,0)[lt]{$\Sigma_{x_1}$}}%
\put(15.4400,-20.5200){\makebox(0,0)[rb]{$M$}}%
\put(24.1900,-20.8000){\makebox(0,0)[lb]{$F_1$ (focal submanifold)@ }}%
\put(23.9400,-29.8700){\makebox(0,0)[lt]{$x_1$}}%
%
\special{pn 8}%
\special{pa 2512 2946}%
\special{pa 2536 2926}%
\special{pa 2560 2904}%
\special{pa 2586 2886}%
\special{pa 2610 2866}%
\special{pa 2636 2846}%
\special{pa 2662 2828}%
\special{pa 2690 2810}%
\special{pa 2716 2794}%
\special{pa 2744 2780}%
\special{pa 2772 2764}%
\special{pa 2802 2748}%
\special{pa 2830 2736}%
\special{pa 2860 2722}%
\special{pa 2890 2710}%
\special{pa 2920 2700}%
\special{pa 2950 2688}%
\special{pa 2980 2678}%
\special{pa 3010 2670}%
\special{pa 3042 2662}%
\special{pa 3072 2654}%
\special{pa 3104 2646}%
\special{pa 3136 2640}%
\special{pa 3166 2636}%
\special{pa 3198 2632}%
\special{pa 3230 2628}%
\special{pa 3262 2624}%
\special{pa 3294 2622}%
\special{pa 3326 2620}%
\special{pa 3358 2620}%
\special{pa 3390 2620}%
\special{pa 3422 2622}%
\special{pa 3454 2624}%
\special{pa 3486 2626}%
\special{pa 3518 2630}%
\special{pa 3550 2634}%
\special{pa 3580 2638}%
\special{pa 3612 2646}%
\special{pa 3644 2652}%
\special{pa 3674 2658}%
\special{pa 3706 2666}%
\special{pa 3736 2676}%
\special{pa 3766 2686}%
\special{pa 3798 2696}%
\special{pa 3828 2706}%
\special{pa 3856 2718}%
\special{pa 3886 2732}%
\special{pa 3916 2744}%
\special{pa 3944 2758}%
\special{pa 3972 2774}%
\special{pa 4000 2790}%
\special{pa 4028 2806}%
\special{pa 4056 2822}%
\special{pa 4082 2840}%
\special{pa 4108 2858}%
\special{pa 4134 2876}%
\special{pa 4160 2896}%
\special{pa 4184 2916}%
\special{pa 4194 2924}%
\special{sp 0.070}%
%
\special{pn 8}%
\special{pa 4166 2920}%
\special{pa 4240 3006}%
\special{fp}%
\special{sh 1}%
\special{pa 4240 3006}%
\special{pa 4212 2942}%
\special{pa 4204 2964}%
\special{pa 4182 2968}%
\special{pa 4240 3006}%
\special{fp}%
\put(31.7700,-25.4000){\makebox(0,0)[lb]{in fact}}%
%
\special{pn 20}%
\special{sh 1}%
\special{ar 4912 3196 10 10 0  6.28318530717959E+0000}%
\special{sh 1}%
\special{ar 4912 3196 10 10 0  6.28318530717959E+0000}%
\put(48.7200,-32.3500){\makebox(0,0)[rb]{{\small $x_1$}}}%
\put(38.1000,-37.5000){\makebox(0,0)[lt]{$\,$}}%
%
\special{pn 8}%
\special{pa 2364 2936}%
\special{pa 2352 2964}%
\special{pa 2334 2990}%
\special{pa 2308 3008}%
\special{pa 2278 3024}%
\special{pa 2248 3034}%
\special{pa 2220 3046}%
\special{pa 2188 3052}%
\special{pa 2156 3058}%
\special{pa 2124 3062}%
\special{pa 2092 3066}%
\special{pa 2060 3066}%
\special{pa 2028 3068}%
\special{pa 1996 3066}%
\special{pa 1964 3064}%
\special{pa 1934 3062}%
\special{pa 1902 3058}%
\special{pa 1870 3052}%
\special{pa 1838 3048}%
\special{pa 1808 3038}%
\special{pa 1776 3032}%
\special{pa 1746 3022}%
\special{pa 1716 3008}%
\special{pa 1688 2994}%
\special{pa 1660 2978}%
\special{pa 1636 2958}%
\special{pa 1616 2932}%
\special{pa 1604 2904}%
\special{pa 1602 2890}%
\special{sp}%
%
\special{pn 8}%
\special{ar 1982 2870 380 142  3.1415927 4.6394417}%
%
\special{pn 8}%
\special{pa 2040 2748}%
\special{pa 2072 2752}%
\special{pa 2104 2756}%
\special{pa 2136 2760}%
\special{pa 2166 2766}%
\special{pa 2198 2774}%
\special{pa 2230 2780}%
\special{pa 2260 2790}%
\special{pa 2290 2802}%
\special{pa 2318 2816}%
\special{pa 2344 2834}%
\special{pa 2366 2856}%
\special{pa 2382 2884}%
\special{pa 2380 2908}%
\special{sp}%
%
\special{pn 8}%
\special{pa 2458 2476}%
\special{pa 2448 2506}%
\special{pa 2430 2532}%
\special{pa 2402 2548}%
\special{pa 2376 2564}%
\special{pa 2346 2576}%
\special{pa 2316 2588}%
\special{pa 2286 2596}%
\special{pa 2254 2600}%
\special{pa 2222 2606}%
\special{pa 2190 2610}%
\special{pa 2158 2612}%
\special{pa 2126 2612}%
\special{pa 2094 2612}%
\special{pa 2062 2610}%
\special{pa 2030 2608}%
\special{pa 1998 2604}%
\special{pa 1968 2598}%
\special{pa 1936 2594}%
\special{pa 1904 2588}%
\special{pa 1874 2578}%
\special{pa 1844 2568}%
\special{pa 1814 2556}%
\special{pa 1784 2544}%
\special{pa 1756 2528}%
\special{pa 1728 2512}%
\special{pa 1704 2492}%
\special{pa 1684 2466}%
\special{pa 1672 2438}%
\special{pa 1670 2414}%
\special{sp}%
%
\special{pn 8}%
\special{pa 1678 2438}%
\special{pa 1688 2408}%
\special{pa 1708 2384}%
\special{pa 1732 2364}%
\special{pa 1760 2348}%
\special{pa 1790 2336}%
\special{pa 1820 2326}%
\special{pa 1850 2318}%
\special{pa 1882 2310}%
\special{pa 1914 2306}%
\special{pa 1946 2302}%
\special{pa 1978 2298}%
\special{pa 2008 2298}%
\special{pa 2040 2296}%
\special{pa 2058 2296}%
\special{sp}%
%
\special{pn 8}%
\special{ar 2106 2448 362 142  4.8655021 6.2831853}%
\special{ar 2106 2448 362 142  0.0000000 0.0870138}%
%
\special{pn 20}%
\special{sh 1}%
\special{ar 2458 2466 10 10 0  6.28318530717959E+0000}%
\special{sh 1}%
\special{ar 2458 2466 10 10 0  6.28318530717959E+0000}%
\put(24.7600,-25.1300){\makebox(0,0)[lt]{$x_2$}}%
%
\special{pn 13}%
\special{ar 1466 2764 2052 344  4.7189472 5.3525288}%
%
\special{pn 20}%
\special{sh 1}%
\special{ar 2048 2448 10 10 0  6.28318530717959E+0000}%
\special{sh 1}%
\special{ar 2048 2448 10 10 0  6.28318530717959E+0000}%
\put(27.1400,-25.6900){\makebox(0,0)[lb]{$\Sigma_{x_2}$}}%
%
\special{pn 8}%
\special{ar 6198 2750 728 854  3.2026028 4.1782381}%
%
\special{pn 8}%
\special{ar 4244 3064 2064 464  4.7186496 5.3546164}%
%
\special{pn 8}%
\special{ar 4592 2400 2066 464  4.7192090 5.3542178}%
%
\special{pn 13}%
\special{ar 4312 2832 2004 468  4.7188386 5.3539596}%
%
\special{pn 13}%
\special{ar 5424 2992 1158 1078  3.5082189 4.4623239}%
%
\special{pn 13}%
\special{ar 4510 2978 912 1188  5.2410850 6.0331494}%
%
\special{pn 20}%
\special{sh 1}%
\special{ar 4940 2256 10 10 0  6.28318530717959E+0000}%
\special{sh 1}%
\special{ar 4940 2256 10 10 0  6.28318530717959E+0000}%
\put(49.0200,-22.0300){\makebox(0,0)[rt]{{\small $x_2$}}}%
%
\special{pn 13}%
\special{ar 4692 3312 248 1250  2.7665784 2.9147938}%
\put(57.5700,-32.2700){\makebox(0,0)[lt]{$\Sigma_{x_1}$}}%
\put(57.7700,-23.0600){\makebox(0,0)[lt]{$\Sigma_{x_2}$}}%
\put(46.0600,-16.3800){\makebox(0,0)[lb]{$F_1$}}%
%
\special{pn 20}%
\special{sh 1}%
\special{ar 4484 2372 10 10 0  6.28318530717959E+0000}%
\special{sh 1}%
\special{ar 4484 2372 10 10 0  6.28318530717959E+0000}%
%
\special{pn 20}%
\special{sh 1}%
\special{ar 4454 3312 10 10 0  6.28318530717959E+0000}%
\special{sh 1}%
\special{ar 4454 3312 10 10 0  6.28318530717959E+0000}%
%
\special{pn 8}%
\special{ar 2382 3904 848 2538  3.2959690 3.4234490}%
%
\special{pn 8}%
\special{ar 2686 3716 352 1570  3.2439234 3.4370335}%
%
\special{pn 13}%
\special{ar 2800 3480 856 1982  3.1274427 3.2783546}%
%
\special{pn 13}%
\special{ar 5102 3190 656 2256  2.9994599 3.3842586}%
%
\special{pn 13}%
\special{ar 4854 3180 420 1692  3.5060163 4.2296180}%
%
\special{pn 8}%
\special{ar 4892 3650 666 884  3.2590437 3.9133287}%
%
\special{pn 8}%
\special{ar 4930 2692 666 856  3.2429779 4.0679690}%
%
\special{pn 8}%
\special{ar 4902 3650 666 884  4.0176507 4.1981455}%
%
\special{pn 8}%
\special{ar 2096 2118 372 132  6.2831853 6.5217539}%
\special{ar 2096 2118 372 132  6.6648950 6.9034636}%
\special{ar 2096 2118 372 132  7.0466048 7.2851734}%
\special{ar 2096 2118 372 132  7.4283145 7.6583418}%
%
\special{pn 8}%
\special{pa 2182 2240}%
\special{pa 2136 2240}%
\special{da 0.070}%
\special{sh 1}%
\special{pa 2136 2240}%
\special{pa 2202 2260}%
\special{pa 2188 2240}%
\special{pa 2202 2220}%
\special{pa 2136 2240}%
\special{fp}%
%
\special{pn 13}%
\special{ar 5674 3226 420 1692  3.5060163 4.2296180}%
%
\special{pn 13}%
\special{pa 5302 3582}%
\special{pa 5298 3550}%
\special{pa 5296 3518}%
\special{pa 5294 3486}%
\special{pa 5290 3454}%
\special{pa 5288 3424}%
\special{pa 5286 3392}%
\special{pa 5284 3360}%
\special{pa 5282 3328}%
\special{pa 5280 3296}%
\special{pa 5278 3264}%
\special{pa 5278 3232}%
\special{pa 5276 3200}%
\special{pa 5276 3168}%
\special{pa 5274 3136}%
\special{pa 5274 3104}%
\special{pa 5272 3072}%
\special{pa 5272 3040}%
\special{pa 5272 3008}%
\special{pa 5270 2976}%
\special{pa 5270 2944}%
\special{pa 5270 2912}%
\special{pa 5270 2880}%
\special{pa 5270 2848}%
\special{pa 5270 2816}%
\special{pa 5270 2784}%
\special{pa 5272 2752}%
\special{pa 5272 2724}%
\special{sp}%
%
\special{pn 13}%
\special{ar 5536 3384 248 1250  2.7686971 2.9147938}%
%
\special{pn 13}%
\special{ar 5852 2898 882 2888  3.3364454 3.6262325}%
%
\special{pn 13}%
\special{ar 5852 2898 892 2690  3.2547708 3.3108553}%
%
\special{pn 13}%
\special{ar 5862 2720 892 1706  2.8077563 3.1744380}%
%
\special{pn 13}%
\special{ar 5852 2898 852 1656  2.7295122 2.8501359}%
\put(50.1900,-15.1900){\makebox(0,0)[lb]{$F_2$}}%
\put(54.2000,-16.7700){\makebox(0,0)[lb]{$F_3$}}%
%
\special{pn 20}%
\special{sh 1}%
\special{ar 5284 3364 10 10 0  6.28318530717959E+0000}%
\special{sh 1}%
\special{ar 5284 3364 10 10 0  6.28318530717959E+0000}%
%
\special{pn 20}%
\special{sh 1}%
\special{ar 4980 2938 10 10 0  6.28318530717959E+0000}%
\special{sh 1}%
\special{ar 4980 2938 10 10 0  6.28318530717959E+0000}%
%
\special{pn 20}%
\special{sh 1}%
\special{ar 5304 2422 10 10 0  6.28318530717959E+0000}%
\special{sh 1}%
\special{ar 5304 2422 10 10 0  6.28318530717959E+0000}%
%
\special{pn 20}%
\special{sh 1}%
\special{ar 5020 2006 10 10 0  6.28318530717959E+0000}%
\special{sh 1}%
\special{ar 5020 2006 10 10 0  6.28318530717959E+0000}%
\end{picture}%
\hspace{0.5truecm}}

\vspace{0.35truecm}

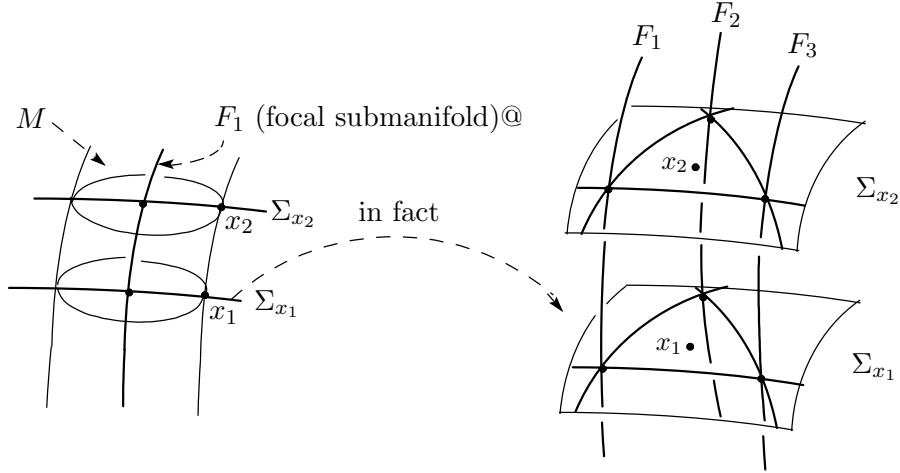
\centerline{{\bf Figure 3$\,:\,$ The set of vertices and focal submanifolds.}}

\vspace{0.5truecm}

We shall show that $(\mathcal B_M,\mathcal O_M)$ is an irreducible topological Tits building of spherical type and satisfies the following 
six conditions:

\vspace{0.25truecm}

(A-I)\ $\displaystyle{\Delta_M=\mathop{\cup}_{x\in M}\mathcal A_x}$, where $\mathcal A_x$ is a coxeter complex defined by 
the simplicial decomposition of the section $\Sigma_x$ of $M$ thorough $x$ by focal walls in $\Sigma_x$;

(A-II)\ 
$|\Delta_M|$ equals to $G/K$, where $|\Delta_M|$ denotes 
$\displaystyle{\mathop{\cup}_{\sigma\in\Delta_M}\sigma}$;

(A-III)\ The set $\mathcal V_M$ of all vertices of $\Delta_M$ equals to the sum of $(r+1)$-pieces of focal submanifolds of $M$ 
and $\mathcal O_M$ is the relative topology (whcih is a metric topology) of focal submanifolds induced from the topology of $G/K$;

(A-IV)\ $\{\mathcal A_x\}_{x\in M}\subset\mathcal A_M$;

(A-V)\ For each $p\in\mathcal V_M$, 
there exists a subfamily $\{\mathcal A_{\lambda}\}_{\lambda\in\Lambda'}$ of $\mathcal A_M$ with $p\in\mathcal A_{\lambda}$ 
($\lambda\in\Lambda'$) such that $\displaystyle{\mathop{\cup}_{\lambda\in\Lambda'}|\mathcal A_{\lambda}|}$ equals to the normal umbrella of 
the focal submanifold of $M$ through $p$ at $p$, where $|\mathcal A_{\lambda}|$ denotes 
$\displaystyle{\mathop{\cup}_{\sigma\in\mathcal A_{\lambda}}\sigma}$;

(A-VI)\ For any points $p,q\in G/K$, there exists a piecewise smooth one-parameter famiy $\{\mathcal A_{\lambda(s)}\}_{s\in[0,1]}$ of $\mathcal A_M$ 
satisfying the following conditions:

($\ast_1$)\ $p\in\mathcal A_{\lambda(0)}$ and $q\in\mathcal A_{\lambda(1)}$;

($\ast_2$)\ There exists a division $0=s_0<s_1<s_2<\cdots<s_k=1$ of $[0,1]$ such that 
$\{\mathcal A_{\lambda(s)}\}_{s\in[s_{i-1},s_i]}$ ($i=1,\cdots,k$) are smooth families and that, for each $i\in\{1,\cdots,k\}$, 
$\displaystyle{\mathop{\cap}_{s\in[s_{i-1},s_i]}|\mathcal A_{]\lambda(s)}|}$ is a common $(r-1)$-simplex of $\mathcal A_{\lambda(s)}$ 
($s\in[s_{i-1},s_i]$), where $r$ is the codimension of $M$.  

\vspace{0.25truecm}

Note that ``smoothness'' of $\{\mathcal A_{\lambda(s)}\}_{s\in[s_{i-1},s_i]}$ in ($\ast_2$) means that the tangent spaces of 
$\mathcal A_{\lambda(s)}$'s ($s\in[s_{i-1},s_i]$) at any common point $p$ give a $C^{\infty}$-curve in the Grassmann manifold $G_r(T_pM)$ 
consisting of $r$-dimensional vector subspaces of $T_pM$.  
Since the apartments $\mathcal A_x$'s ($x\in M$) are coxeter complexes, all apartments $\mathcal A_{\lambda}$ ($\lambda\in\Lambda$) also are 
coxeter complexes.  Hence the condition (B2) holds.  

Let $\{\mathcal H^x_1,\cdots,\mathcal H^x_k\}$ be the set of all focal walls in $\Sigma_x$.  
Take any $i\in\{1,\cdots,k\}$.  Take a parallel normal vector field $\widetilde{\vv}_i$ of $M$ with 
$\displaystyle{\exp^{\perp}((\widetilde{\vv}_i)_x)\in\mathcal H^x_i\setminus\left(\mathop{\cup}_{j\in\{1,\cdots,k\}\setminus\{i\}}\mathcal H^x_j
\right)}$.  For any $y\in L^{\widetilde{\vv}_i}_x$, $\Sigma_x\cap\Sigma_y=\mathcal H^x_i$ holds.  
From this fact, it follows that each $(r-1)$-dimensional simplex of $\Delta_M$ is contained in infinitely many chambers.  
Hence the condition (B1) holds.  

Take any two simplices $\sigma_1$ and $\sigma_2$ of $\Delta_M$.  
Let $x_i$ ($i=1,2$) be points of $M$ such that $\sigma_i$ is contained in the chamber $C_i$ of $\mathcal A_{x_i}$ containing $x_i$.  
Let $\widetilde{\vv}_i$ ($i=1,\cdots,k$) be the above focal normal vector field of $M$.  
Denote by $\mathfrak F_i$ the foliation on $M$ given by the fibres of the focal map $\eta_{\widetilde{\vv}_i}:M\to F_{\widetilde{\vv}_i}$ 
(which is a submersion).  The family $\{\mathfrak F_i\}_{i=1}^k$ of these foliations gives a net on $M$.  
Here ``net'' means that $\displaystyle{TM=\mathop{\oplus}_{i=1}^kT\mathcal D_{\mathfrak F_i}}$ holds (the terminology ``net'' was originally 
used in \cite{Ko1} (in \cite{RS} and \cite{Ko2} later), where $\mathcal D_{\mathfrak F_i}$ is the distribution on $M$ associated to 
$\mathfrak F_i$ (i.e., $(\mathcal D_{\mathfrak F_i})_x=T_xL^{\widetilde{\vv}_i}_x$ ($x\in M$)).  
Since $M$ is compact, it is shown that there exists a piecewise smooth curve $\alpha:[0,1]\to M$ with $\alpha(0)=x_1$ and $\alpha(1)=x_2$ 
which admits a division $0=s_0<s_1<\cdots<s_l=1$ of $[0,1]$ satisfying 

\vspace{0.15truecm}

$(\ast)$\ $\alpha|_{[s_{i-1},s_i]}$ ($i=1,\cdots,l$) are $C^{\infty}$-curves in a leaf of some $\mathfrak F_{j(i)}$.  

\vspace{0.15truecm}

\noindent
Set $\Sigma_s:=\Sigma_{\alpha(s)}$ ($s\in[0,1]$).  Then we can show that $\{\Sigma_s\}_{s\in[s_{i-1},s_i]}$ share the focal wall 
$\mathcal H^{\alpha(s_{i-1})}_{j(i)}$ ($i=1,\cdots,l$) (see Figure 4).  From this fact, we can find an apartment of $\mathcal B_M$ containing 
both $\sigma_1$ and $\sigma_2$.  
Hence the condition (B3) holds.  
Also, by this discussion, it is shown that $(\mathcal B_M,\mathcal O_M)$ satisfies the condition (A-VI).  

According to the construction of $\mathcal A$, it is clear that the conditions (B4) and (B5) hold.  
Also, since $\mathcal V_M$ consists of $(r+1)$-pieces of focal submanifolds of $M$, it follows from the definition of $\mathcal S_M$ that 
the condition (TB2) holds.  Therefore $(\mathcal B_M,\mathcal O_M)$ is a topological Tits building of spherical type.  
It is clear that $(\mathcal B_M,\mathcal O_M)$ is of rank $r$ and irreducible.  
Also, it is clear that this topological Tits building $(\mathcal B_M,\mathcal O_M)$ satisfies the conditions (A-I)$-$(A-IV).  
We shall show that $(\mathcal B_M,\mathcal O_M)$ satisfies the conditions (A-V).  
Take any vertex $p$ of $(\mathcal B_M,\mathcal O_M)$.  
Let $x$ be a point of $M$ such that $p$ is a vertex of the chamber $C$ of $\mathcal A_x$ containing $x$, and 
$\widetilde{\vv}$ be a focal normal vector field of $M$ satisfying $\eta_{\widetilde{\vv}}(x)=p$.  
According to the proof of Theorem C (see Section 4), the normal umbrella $S_p:=\exp^{\perp}(T^{\perp}_pF_{\widetilde{\vv}})$ of 
the focal submanifold $F_{\widetilde{\vv}}:=\eta_{\widetilde{\vv}}(M)$ equals to 
$\displaystyle{\mathop{\cup}_{y\in L^{\widetilde{\vv}}_x}\Sigma_y}$, where $L^{\widetilde{\vv}}_x:=\eta_{\widetilde{\vv}}^{-1}(p)$.  
Hence we obtain $\displaystyle{\mathop{\cup}_{y\in L^{\widetilde{\vv}}_x}|\mathcal A_y|=S_p}$.  
Therefore $\{\mathcal A_y\}_{y\in L^{\widetilde{\vv}}_x}$ is a subfamily of $\mathcal A_M$ as in (A-V).  
Thus $(\mathcal B_M,\mathcal O_M)$ satisfies the conditions (A-V).  Thus $(\mathcal B_M,\mathcal O_M)$ satisfies the above conditions (A-I)$-$(A-VI).  
We call $(\mathcal B_M,\mathcal O_M)$ the {\it topological Tits building of spherical type associated to} $M$.  

\vspace{0.5truecm}

\centerline{
\unitlength 0.1in
\begin{picture}( 26.1300, 21.8100)( 15.0000,-30.5700)
%
\special{pn 8}%
\special{ar 2628 2754 756 794  3.2035119 4.1797302}%
%
\special{pn 8}%
\special{ar 4194 2830 754 796  3.2026964 4.1775947}%
%
\special{pn 8}%
\special{pa 1876 2708}%
\special{pa 1908 2708}%
\special{pa 1940 2708}%
\special{pa 1972 2706}%
\special{pa 2004 2704}%
\special{pa 2036 2702}%
\special{pa 2068 2702}%
\special{pa 2100 2702}%
\special{pa 2132 2702}%
\special{pa 2164 2700}%
\special{pa 2196 2700}%
\special{pa 2228 2700}%
\special{pa 2260 2700}%
\special{pa 2292 2700}%
\special{pa 2324 2698}%
\special{pa 2356 2698}%
\special{pa 2388 2698}%
\special{pa 2420 2700}%
\special{pa 2452 2700}%
\special{pa 2484 2700}%
\special{pa 2516 2702}%
\special{pa 2548 2702}%
\special{pa 2580 2704}%
\special{pa 2612 2704}%
\special{pa 2644 2704}%
\special{pa 2676 2706}%
\special{pa 2708 2708}%
\special{pa 2740 2708}%
\special{pa 2772 2708}%
\special{pa 2804 2712}%
\special{pa 2836 2714}%
\special{pa 2868 2716}%
\special{pa 2900 2718}%
\special{pa 2932 2720}%
\special{pa 2964 2722}%
\special{pa 2996 2724}%
\special{pa 3028 2728}%
\special{pa 3060 2730}%
\special{pa 3092 2732}%
\special{pa 3124 2736}%
\special{pa 3154 2740}%
\special{pa 3186 2742}%
\special{pa 3218 2746}%
\special{pa 3250 2750}%
\special{pa 3282 2752}%
\special{pa 3314 2756}%
\special{pa 3346 2760}%
\special{pa 3378 2764}%
\special{pa 3410 2768}%
\special{pa 3440 2774}%
\special{pa 3442 2774}%
\special{sp}%
%
\special{pn 13}%
\special{ar 1916 3132 2554 646  4.7188257 5.3539962}%
%
\special{pn 13}%
\special{ar 3390 3052 1198 1004  3.5080270 4.4622196}%
%
\special{pn 13}%
\special{ar 2446 3044 946 1108  5.2403136 6.0323492}%
%
\special{pn 8}%
\special{ar 2430 3066 156 1668  3.3649418 4.3360191}%
%
\special{pn 8}%
\special{pa 2360 1542}%
\special{pa 2382 1520}%
\special{pa 2404 1496}%
\special{pa 2428 1476}%
\special{pa 2452 1452}%
\special{pa 2476 1434}%
\special{pa 2502 1414}%
\special{pa 2528 1394}%
\special{pa 2554 1376}%
\special{pa 2580 1358}%
\special{pa 2608 1340}%
\special{pa 2634 1324}%
\special{pa 2664 1310}%
\special{pa 2692 1294}%
\special{pa 2720 1280}%
\special{pa 2748 1266}%
\special{pa 2778 1252}%
\special{pa 2806 1238}%
\special{pa 2836 1228}%
\special{pa 2866 1216}%
\special{pa 2896 1206}%
\special{pa 2926 1194}%
\special{pa 2956 1184}%
\special{pa 2986 1174}%
\special{pa 3018 1166}%
\special{pa 3048 1156}%
\special{pa 3080 1148}%
\special{pa 3110 1142}%
\special{pa 3142 1132}%
\special{pa 3174 1128}%
\special{pa 3204 1120}%
\special{pa 3236 1114}%
\special{pa 3268 1108}%
\special{pa 3282 1106}%
\special{sp}%
%
\special{pn 13}%
\special{sh 1}%
\special{ar 2922 2400 10 10 0  6.28318530717959E+0000}%
\special{sh 1}%
\special{ar 2922 2400 10 10 0  6.28318530717959E+0000}%
%
\special{pn 20}%
\special{sh 1}%
\special{ar 2668 1974 10 10 0  6.28318530717959E+0000}%
\special{sh 1}%
\special{ar 2668 1974 10 10 0  6.28318530717959E+0000}%
%
\special{pn 8}%
\special{pa 3706 1010}%
\special{pa 3104 1220}%
\special{da 0.070}%
\special{sh 1}%
\special{pa 3104 1220}%
\special{pa 3174 1218}%
\special{pa 3154 1202}%
\special{pa 3160 1180}%
\special{pa 3104 1220}%
\special{fp}%
\put(40.8100,-21.4800){\makebox(0,0)[lt]{$\Sigma_{s_{i-1}}$}}%
\put(37.4500,-10.4600){\makebox(0,0)[lb]{$\Sigma_{s_i}$}}%
\put(15.1000,-14.3200){\makebox(0,0)[lb]{$L^{\widetilde{\vv}_{j(i)}}_{\alpha(s_{i-1})}$}}%
\put(37.1500,-16.8000){\makebox(0,0)[lb]{$\alpha(s_i)$}}%
\put(39.3300,-24.6000){\makebox(0,0)[lt]{$\alpha(s_{i-1})$}}%
%
\special{pn 8}%
\special{pa 2242 2084}%
\special{pa 2274 2082}%
\special{pa 2306 2082}%
\special{pa 2338 2082}%
\special{pa 2370 2080}%
\special{pa 2402 2078}%
\special{pa 2434 2078}%
\special{pa 2466 2078}%
\special{pa 2498 2078}%
\special{pa 2530 2076}%
\special{pa 2562 2076}%
\special{pa 2594 2076}%
\special{pa 2626 2076}%
\special{pa 2658 2076}%
\special{pa 2690 2074}%
\special{pa 2722 2074}%
\special{pa 2754 2074}%
\special{pa 2786 2076}%
\special{pa 2818 2076}%
\special{pa 2850 2076}%
\special{pa 2882 2078}%
\special{pa 2914 2078}%
\special{pa 2946 2078}%
\special{pa 2978 2078}%
\special{pa 3010 2080}%
\special{pa 3042 2082}%
\special{pa 3074 2082}%
\special{pa 3106 2084}%
\special{pa 3138 2086}%
\special{pa 3170 2088}%
\special{pa 3202 2090}%
\special{pa 3234 2092}%
\special{pa 3266 2092}%
\special{pa 3298 2096}%
\special{pa 3330 2098}%
\special{pa 3362 2100}%
\special{pa 3394 2102}%
\special{pa 3426 2106}%
\special{pa 3458 2110}%
\special{pa 3490 2112}%
\special{pa 3520 2114}%
\special{pa 3552 2118}%
\special{pa 3584 2122}%
\special{pa 3616 2124}%
\special{pa 3648 2128}%
\special{pa 3680 2132}%
\special{pa 3712 2136}%
\special{pa 3744 2142}%
\special{pa 3776 2144}%
\special{pa 3806 2148}%
\special{pa 3808 2148}%
\special{sp}%
%
\special{pn 13}%
\special{ar 2696 2396 228 424  4.7123890 6.2831853}%
%
\special{pn 13}%
\special{ar 2696 2406 298 432  3.7233884 4.7371364}%
%
\special{pn 20}%
\special{sh 1}%
\special{ar 2450 2168 10 10 0  6.28318530717959E+0000}%
\special{sh 1}%
\special{ar 2450 2168 10 10 0  6.28318530717959E+0000}%
%
\special{pn 20}%
\special{sh 1}%
\special{ar 2658 2268 10 10 0  6.28318530717959E+0000}%
\special{sh 1}%
\special{ar 2658 2268 10 10 0  6.28318530717959E+0000}%
\put(19.9400,-29.8400){\makebox(0,0)[lt]{$\displaystyle{\mathop{\cap}_{s\in[s_{i-1},s_i]}\Sigma_s=\mathcal H^{\alpha(s_{i-1})}_{j(i)}}$}}%
%
\special{pn 8}%
\special{pa 3714 2072}%
\special{pa 3548 2194}%
\special{da 0.070}%
\special{sh 1}%
\special{pa 3548 2194}%
\special{pa 3614 2172}%
\special{pa 3590 2164}%
\special{pa 3590 2140}%
\special{pa 3548 2194}%
\special{fp}%
%
\special{pn 8}%
\special{ar 3904 2158 222 148  3.6622585 3.9865828}%
\special{ar 3904 2158 222 148  4.1811774 4.5055018}%
\special{ar 3904 2158 222 148  4.7000964 5.0244207}%
\special{ar 3904 2158 222 148  5.2190153 5.5433396}%
\special{ar 3904 2158 222 148  5.7379342 5.9507002}%
%
\special{pn 8}%
\special{pa 1946 1478}%
\special{pa 2558 2002}%
\special{da 0.070}%
\special{sh 1}%
\special{pa 2558 2002}%
\special{pa 2520 1944}%
\special{pa 2518 1966}%
\special{pa 2494 1974}%
\special{pa 2558 2002}%
\special{fp}%
%
\special{pn 20}%
\special{sh 1}%
\special{ar 2746 1984 10 10 0  6.28318530717959E+0000}%
\special{sh 1}%
\special{ar 2746 1984 10 10 0  6.28318530717959E+0000}%
%
\special{pn 8}%
\special{pa 2272 2650}%
\special{pa 2280 2620}%
\special{pa 2288 2588}%
\special{pa 2296 2558}%
\special{pa 2304 2526}%
\special{pa 2312 2496}%
\special{pa 2320 2464}%
\special{pa 2328 2434}%
\special{pa 2338 2404}%
\special{pa 2344 2372}%
\special{pa 2352 2342}%
\special{pa 2362 2310}%
\special{pa 2370 2280}%
\special{pa 2380 2248}%
\special{pa 2388 2218}%
\special{pa 2396 2188}%
\special{pa 2406 2156}%
\special{pa 2416 2126}%
\special{pa 2424 2096}%
\special{pa 2434 2064}%
\special{pa 2442 2034}%
\special{pa 2452 2004}%
\special{pa 2462 1974}%
\special{pa 2472 1942}%
\special{pa 2480 1912}%
\special{pa 2490 1882}%
\special{pa 2502 1852}%
\special{pa 2512 1822}%
\special{pa 2520 1790}%
\special{pa 2530 1760}%
\special{pa 2542 1730}%
\special{pa 2554 1700}%
\special{pa 2566 1670}%
\special{pa 2578 1640}%
\special{pa 2590 1612}%
\special{pa 2602 1582}%
\special{pa 2614 1552}%
\special{pa 2626 1522}%
\special{pa 2634 1508}%
\special{sp}%
%
\special{pn 8}%
\special{pa 2638 1514}%
\special{pa 2662 1494}%
\special{pa 2686 1472}%
\special{pa 2710 1454}%
\special{pa 2736 1434}%
\special{pa 2764 1416}%
\special{pa 2790 1400}%
\special{pa 2818 1384}%
\special{pa 2848 1370}%
\special{pa 2876 1356}%
\special{pa 2906 1344}%
\special{pa 2936 1334}%
\special{pa 2966 1322}%
\special{pa 2996 1312}%
\special{pa 3026 1302}%
\special{pa 3056 1294}%
\special{pa 3088 1286}%
\special{pa 3120 1282}%
\special{pa 3150 1276}%
\special{pa 3182 1270}%
\special{pa 3214 1266}%
\special{pa 3246 1264}%
\special{pa 3278 1262}%
\special{pa 3310 1260}%
\special{pa 3342 1258}%
\special{pa 3374 1258}%
\special{pa 3406 1258}%
\special{pa 3430 1260}%
\special{sp}%
%
\special{pn 8}%
\special{pa 3102 2084}%
\special{pa 3108 2054}%
\special{pa 3116 2022}%
\special{pa 3122 1990}%
\special{pa 3128 1960}%
\special{pa 3136 1928}%
\special{pa 3144 1898}%
\special{pa 3152 1866}%
\special{pa 3160 1836}%
\special{pa 3168 1804}%
\special{pa 3178 1774}%
\special{pa 3186 1744}%
\special{pa 3198 1714}%
\special{pa 3208 1684}%
\special{pa 3218 1652}%
\special{pa 3230 1622}%
\special{pa 3240 1592}%
\special{pa 3252 1562}%
\special{pa 3264 1534}%
\special{pa 3278 1504}%
\special{pa 3294 1476}%
\special{pa 3308 1448}%
\special{pa 3322 1420}%
\special{pa 3338 1392}%
\special{pa 3356 1366}%
\special{pa 3372 1338}%
\special{pa 3392 1312}%
\special{pa 3414 1288}%
\special{pa 3428 1272}%
\special{sp}%
%
\special{pn 8}%
\special{pa 2894 1790}%
\special{pa 2678 1956}%
\special{da 0.070}%
\special{sh 1}%
\special{pa 2678 1956}%
\special{pa 2742 1932}%
\special{pa 2720 1924}%
\special{pa 2718 1900}%
\special{pa 2678 1956}%
\special{fp}%
%
\special{pn 8}%
\special{ar 3676 1956 882 358  3.5942718 3.6911241}%
\special{ar 3676 1956 882 358  3.7492355 3.8460878}%
\special{ar 3676 1956 882 358  3.9041991 4.0010514}%
\special{ar 3676 1956 882 358  4.0591628 4.1560151}%
\special{ar 3676 1956 882 358  4.2141265 4.3109788}%
\special{ar 3676 1956 882 358  4.3690902 4.4659425}%
\special{ar 3676 1956 882 358  4.5240539 4.6209062}%
\special{ar 3676 1956 882 358  4.6790175 4.7123890}%
%
\special{pn 8}%
\special{pa 3706 1286}%
\special{pa 3230 1450}%
\special{da 0.070}%
\special{sh 1}%
\special{pa 3230 1450}%
\special{pa 3300 1448}%
\special{pa 3280 1432}%
\special{pa 3286 1410}%
\special{pa 3230 1450}%
\special{fp}%
\put(37.2500,-13.4900){\makebox(0,0)[lb]{$\Sigma_s$}}%
%
\special{pn 8}%
\special{pa 2944 1908}%
\special{pa 2756 1980}%
\special{da 0.070}%
\special{sh 1}%
\special{pa 2756 1980}%
\special{pa 2826 1974}%
\special{pa 2806 1960}%
\special{pa 2812 1936}%
\special{pa 2756 1980}%
\special{fp}%
%
\special{pn 8}%
\special{ar 3608 1984 758 156  3.5956712 3.7271060}%
\special{ar 3608 1984 758 156  3.8059669 3.9374017}%
\special{ar 3608 1984 758 156  4.0162626 4.1476975}%
\special{ar 3608 1984 758 156  4.2265584 4.3579932}%
\special{ar 3608 1984 758 156  4.4368541 4.5682889}%
\special{ar 3608 1984 758 156  4.6471498 4.7123890}%
\put(36.6600,-18.9100){\makebox(0,0)[lb]{$\alpha(s)$}}%
%
\special{pn 8}%
\special{pa 3072 2324}%
\special{pa 2934 2388}%
\special{da 0.070}%
\special{sh 1}%
\special{pa 2934 2388}%
\special{pa 3004 2378}%
\special{pa 2982 2366}%
\special{pa 2986 2342}%
\special{pa 2934 2388}%
\special{fp}%
%
\special{pn 8}%
\special{ar 3458 2580 614 330  4.0255713 4.1528248}%
\special{ar 3458 2580 614 330  4.2291768 4.3564303}%
\special{ar 3458 2580 614 330  4.4327824 4.5600358}%
\special{ar 3458 2580 614 330  4.6363879 4.7636413}%
\special{ar 3458 2580 614 330  4.8399934 4.9672468}%
\special{ar 3458 2580 614 330  5.0435989 5.1708523}%
\special{ar 3458 2580 614 330  5.2472044 5.3744579}%
\special{ar 3458 2580 614 330  5.4508099 5.5780634}%
\special{ar 3458 2580 614 330  5.6544154 5.7816689}%
%
\special{pn 8}%
\special{pa 2658 2268}%
\special{pa 2666 2238}%
\special{pa 2674 2206}%
\special{pa 2684 2176}%
\special{pa 2694 2146}%
\special{pa 2702 2114}%
\special{pa 2712 2084}%
\special{pa 2720 2054}%
\special{pa 2730 2024}%
\special{pa 2740 1992}%
\special{pa 2752 1962}%
\special{pa 2762 1932}%
\special{pa 2772 1902}%
\special{pa 2784 1872}%
\special{pa 2794 1842}%
\special{pa 2806 1812}%
\special{pa 2818 1782}%
\special{pa 2828 1752}%
\special{pa 2840 1722}%
\special{pa 2854 1694}%
\special{pa 2866 1664}%
\special{pa 2876 1634}%
\special{pa 2890 1604}%
\special{pa 2904 1576}%
\special{pa 2916 1546}%
\special{pa 2930 1518}%
\special{pa 2944 1488}%
\special{pa 2958 1460}%
\special{pa 2972 1432}%
\special{pa 2986 1402}%
\special{pa 3000 1374}%
\special{pa 3016 1346}%
\special{pa 3032 1318}%
\special{pa 3048 1290}%
\special{pa 3050 1286}%
\special{sp}%
%
\special{pn 8}%
\special{pa 2232 2094}%
\special{pa 2264 2102}%
\special{pa 2294 2112}%
\special{pa 2324 2124}%
\special{pa 2354 2136}%
\special{pa 2384 2144}%
\special{pa 2414 2154}%
\special{pa 2444 2166}%
\special{pa 2474 2178}%
\special{pa 2504 2190}%
\special{pa 2534 2202}%
\special{pa 2564 2214}%
\special{pa 2592 2226}%
\special{pa 2622 2238}%
\special{pa 2652 2252}%
\special{pa 2680 2266}%
\special{pa 2708 2280}%
\special{pa 2738 2294}%
\special{pa 2768 2306}%
\special{pa 2796 2320}%
\special{pa 2824 2336}%
\special{pa 2852 2350}%
\special{pa 2880 2366}%
\special{pa 2908 2382}%
\special{pa 2936 2398}%
\special{pa 2964 2414}%
\special{pa 2990 2432}%
\special{pa 3018 2448}%
\special{pa 3044 2466}%
\special{pa 3072 2482}%
\special{pa 3098 2502}%
\special{pa 3124 2520}%
\special{pa 3150 2538}%
\special{pa 3174 2558}%
\special{pa 3200 2578}%
\special{pa 3224 2598}%
\special{pa 3248 2618}%
\special{pa 3274 2638}%
\special{pa 3296 2662}%
\special{pa 3320 2684}%
\special{pa 3342 2706}%
\special{pa 3364 2730}%
\special{pa 3386 2754}%
\special{pa 3406 2778}%
\special{pa 3410 2784}%
\special{sp}%
%
\special{pn 20}%
\special{sh 1}%
\special{ar 3370 2506 10 10 0  6.28318530717959E+0000}%
\special{sh 1}%
\special{ar 3350 2388 10 10 0  6.28318530717959E+0000}%
\special{sh 1}%
\special{ar 3330 2286 10 10 0  6.28318530717959E+0000}%
\special{sh 1}%
\special{ar 3300 2176 10 10 0  6.28318530717959E+0000}%
\special{sh 1}%
\special{ar 3250 2030 10 10 0  6.28318530717959E+0000}%
\special{sh 1}%
\special{ar 3222 1920 10 10 0  6.28318530717959E+0000}%
\special{sh 1}%
\special{ar 3152 1744 10 10 0  6.28318530717959E+0000}%
\special{sh 1}%
\special{ar 3152 1744 10 10 0  6.28318530717959E+0000}%
\special{sh 1}%
\special{ar 3082 1616 10 10 0  6.28318530717959E+0000}%
\special{sh 1}%
\special{ar 3024 1534 10 10 0  6.28318530717959E+0000}%
\special{sh 1}%
\special{ar 3024 1534 10 10 0  6.28318530717959E+0000}%
%
\special{pn 8}%
\special{ar 4130 2278 1474 2618  3.1476733 3.4462538}%
%
\special{pn 8}%
\special{ar 4130 2278 1474 2526  3.4909290 3.5584670}%
%
\special{pn 8}%
\special{ar 4328 2084 1226 1956  3.1415927 3.5459993}%
%
\special{pn 8}%
\special{pa 3222 1248}%
\special{pa 3236 1220}%
\special{pa 3248 1192}%
\special{pa 3262 1164}%
\special{pa 3278 1136}%
\special{pa 3290 1112}%
\special{sp}%
%
\special{pn 8}%
\special{ar 1956 2572 456 486  1.5707963 1.6983202}%
\special{ar 1956 2572 456 486  1.7748346 1.9023585}%
\special{ar 1956 2572 456 486  1.9788728 2.1063968}%
\special{ar 1956 2572 456 486  2.1829111 2.3104350}%
\special{ar 1956 2572 456 486  2.3869494 2.5144733}%
\special{ar 1956 2572 456 486  2.5909876 2.7185115}%
\special{ar 1956 2572 456 486  2.7950259 2.9225498}%
\special{ar 1956 2572 456 486  2.9990641 3.1265880}%
\special{ar 1956 2572 456 486  3.2031024 3.3306263}%
\special{ar 1956 2572 456 486  3.4071406 3.5346646}%
\special{ar 1956 2572 456 486  3.6111789 3.7387028}%
\special{ar 1956 2572 456 486  3.8152172 3.9427411}%
\special{ar 1956 2572 456 486  4.0192554 4.1467793}%
\special{ar 1956 2572 456 486  4.2232937 4.3508176}%
\special{ar 1956 2572 456 486  4.4273319 4.5548558}%
\special{ar 1956 2572 456 486  4.6313702 4.7588941}%
\special{ar 1956 2572 456 486  4.8354084 4.9629324}%
\special{ar 1956 2572 456 486  5.0394467 5.1669706}%
\special{ar 1956 2572 456 486  5.2434850 5.3710089}%
%
\special{pn 8}%
\special{pa 2262 2212}%
\special{pa 2480 2388}%
\special{da 0.070}%
\special{sh 1}%
\special{pa 2480 2388}%
\special{pa 2440 2330}%
\special{pa 2438 2354}%
\special{pa 2414 2362}%
\special{pa 2480 2388}%
\special{fp}%
%
\special{pn 8}%
\special{ar 3448 2864 486 488  4.0624699 4.1859267}%
\special{ar 3448 2864 486 488  4.2600007 4.3834575}%
\special{ar 3448 2864 486 488  4.4575316 4.5809884}%
\special{ar 3448 2864 486 488  4.6550625 4.7785192}%
\special{ar 3448 2864 486 488  4.8525933 4.9760501}%
\special{ar 3448 2864 486 488  5.0501242 5.1735810}%
\special{ar 3448 2864 486 488  5.2476551 5.3711118}%
\special{ar 3448 2864 486 488  5.4451859 5.5686427}%
\special{ar 3448 2864 486 488  5.6427168 5.7661736}%
\special{ar 3448 2864 486 488  5.8402476 5.9240933}%
%
\special{pn 8}%
\special{pa 3162 2480}%
\special{pa 3122 2516}%
\special{da 0.070}%
\special{sh 1}%
\special{pa 3122 2516}%
\special{pa 3184 2486}%
\special{pa 3162 2480}%
\special{pa 3158 2456}%
\special{pa 3122 2516}%
\special{fp}%
\put(38.6300,-27.3600){\makebox(0,0)[lt]{$\gamma_{\widetilde{\vv}_{j(i)}}$}}%
%
\special{pn 4}%
\special{pa 2450 2498}%
\special{pa 2426 2474}%
\special{dt 0.027}%
\special{pa 2510 2498}%
\special{pa 2450 2444}%
\special{dt 0.027}%
\special{pa 2576 2504}%
\special{pa 2484 2418}%
\special{dt 0.027}%
\special{pa 2640 2510}%
\special{pa 2514 2392}%
\special{dt 0.027}%
\special{pa 2704 2514}%
\special{pa 2544 2364}%
\special{dt 0.027}%
\special{pa 2770 2518}%
\special{pa 2578 2340}%
\special{dt 0.027}%
\special{pa 2834 2524}%
\special{pa 2612 2318}%
\special{dt 0.027}%
\special{pa 2902 2532}%
\special{pa 2648 2296}%
\special{dt 0.027}%
\special{pa 2974 2544}%
\special{pa 2684 2274}%
\special{dt 0.027}%
\special{pa 3034 2544}%
\special{pa 2718 2252}%
\special{dt 0.027}%
\special{pa 3110 2560}%
\special{pa 2754 2230}%
\special{dt 0.027}%
\special{pa 3178 2568}%
\special{pa 2790 2208}%
\special{dt 0.027}%
\special{pa 3246 2576}%
\special{pa 2830 2190}%
\special{dt 0.027}%
\special{pa 3294 2564}%
\special{pa 2868 2168}%
\special{dt 0.027}%
\special{pa 3240 2460}%
\special{pa 2904 2148}%
\special{dt 0.027}%
\special{pa 3138 2310}%
\special{pa 2956 2142}%
\special{dt 0.027}%
\end{picture}%
\hspace{1truecm}}

\vspace{0.5truecm}

\centerline{{\bf Figure 4$\,:\,$ The behaviour of $\Sigma_s$ ($s\in[s_{i-1},s_i]$)}}

\vspace{0.5truecm}

Now we shall prove Theorem A in the case where the sections are spheres by using the associated topological Tits building 
$(\mathcal B_M,\mathcal O_M)$.  

\vspace{0.5truecm}

\noindent
{\it Proof of Theorem A (spherical section-case).}\ \ 
Suppose that the codimension of $M$ is greater than two.  Denote by $r$ the codimension of $M$.  
Since the topological Tits buildings $(\mathcal B_M,\mathcal O_M)$ of spherical type associated to $M$ is irreducible and of rank greater than two, 
it follows from the classification theorem of irreducible topological Tits buildings in \cite{BS} that $(\mathcal B_M,\mathcal O_M)$ is isomorphic to 
the topological Tits building $(\mathcal B_{G'},\mathcal O_{G'})$ associated to a simple non-compact Lie group $G'$ without the center, which is 
defined as a toplological building of spherical type having parabolic subgroups of $G$ as simplices.  
On the other hand, according to the result in \cite{Th} stated in Introduction, $(\mathcal B_{G'},\mathcal O_{G'})$ is isomorphic 
to the topological Tits building $(\mathcal B_{K'\cdot p},\mathcal O_{K'\cdot p})$ of spherical type associated to a principal orbit 
$K'\cdot p$ (which is a full irreducible isoparametric submanifold of codimension $(r+1)$) of the isotropy representation 
of the irreducible symmetric space $G'/K'$ of non-compact type and rank greater $(r+1)$, where $K'$ is the maximal compact Lie subgroup of $G'$.  
Hence $(\mathcal B_M,\mathcal O_M)$ is isomorphic to $(\mathcal B_{K'\cdot p},\mathcal O_{K'\cdot p})$ as topological Tits building.  
Let $\psi$ be an isomorphism of $(\mathcal B_M,\mathcal O_M)$ onto $(\mathcal B_{K'\cdot p},\mathcal O_{K'\cdot p})$.  
For each chamber $C\in(\mathcal S_M)_r$, there uniquely exists the homothety $\widetilde{\psi}_C$ of $|C|$ onto $|\psi(C)|$ such that 
$\widetilde{\psi}_C(|\sigma|)=|\psi(\sigma)|$ holds for all $\sigma\in\mathcal S_M$ with $|\sigma|\subset|C|$, where we note that 
both $C$ and $\psi(C)$ are chambers of the same kind of coxeter complex defined on $r$-dimensional spheres.  Furthermore, by patching 
$\widetilde{\psi}_C$'s ($C\in(\mathcal S_M)_r$), we can construct a map $\widetilde{\psi}$ of $G/K$ onto the unit sphere $S^{n+r}(1)$.  
It is easy to show that $\widetilde{\psi}$ is a bijection of $G/K$ onto $S^{n+r}(1)$ and that $\widetilde{\psi}|_{\Sigma_x}$ is a homothety of 
$\Sigma_x$ onto a totally geodesic sphere (which is the base space of an apartment of $(\mathcal B_{K'\cdot p},\mathcal O_{K'\cdot p})$) 
for all $x\in M$.  Let $F_i$ ($i=1,\cdots,r+1$) be focal submanifolds of $M$ giving $\mathcal B_M$ and $F'_i$ ($i=1,\cdots,r+1$) 
focal submanifolds of $K'\cdot p$ giving $\mathcal B_{K'\cdot p}$.  
Denote by $\mathfrak F$ (resp. $\widehat{\mathfrak F}$) the singular Riemannian foliation consisting of parallel submanifolds and focal submanifolds 
of $M$ (resp. $K'\cdot p$).  Since $\psi$ is an isomorphism between the topological Tits buildings $(\mathcal B_M,\mathcal O_M)$ and 
$(\mathcal B_{K'\cdot p},\mathcal O_{K'\cdot p})$, it follows from the construction of $\widetilde{\psi}$ 
that $\widetilde{\psi}|_{F_i}$ ($i=1,\cdots,r+1$) is a homeomorphism of $F_i$ onto some $F'_{j(i)}$ and that 
$\widetilde{\psi}$ maps the leaves of $F$ to the leaves of $\mathfrak F'$.  From the above facts, we can derive that 
$\widetilde{\psi}:G/K\to S^{n+r}(1)$ is a homeomorphism.  
Furthermore, it follows from this fact that $G/K$ is homothetic to $S^{n+r}(1)$.  This contradicts that $G/K$ is of rank greater than one.  
Therefore we obtain $r=2$.  \qed

\section{A new proof of Theorem A (projective section-case)} 
In this section, we give a new proof of Theorem A in the case where the sections of the equifocal submanifold are real projective spaces.  
First, for an equifocal submanifold in a symmetric space of compact type, we introduce a new topology of the symmetric space.  
Let $M$ be an equifocal submanifold in a symmetric space $G/K$ of compact type and $\mathfrak F$ the polar foliation cosisting of the parallel submanifold 
and the focal submanifold of $M$.  For each point $p$ of $G/K$, denote by $M_p$ the parallel submanifold (or the focal submanifold) of $M$ through $p$.  
The normal umbrella $S_p:=\exp^{\perp}(T_p^{\perp}M_p)$ is the {\it slice} of the polar foliation $\mathfrak F$ {\it at} $p$.  
Denote by $\mathcal O_p$ the topology of $S_p$.  Let $\mathcal O_{\rm sl}$ be the topology of $G/K$ generated by 
$\mathop{\cup}_{p\in G/K}\{U\in\mathcal O_p\,|\,p\in U\}$.  We call this topology $\mathcal O_{\rm sl}$ the {\it slice topology of} $G/K$ 
{\it associated to} $M$.  It is easy to show that $\mathcal O_{\rm sl}$ is equal to the topology generated by 
$\displaystyle{\mathop{\cup}_{p\in M}\mathcal O_p}$ (see Figure 5), where $\mathcal O_p$ is the topology of the section $\Sigma_p$ of $M$ through $p$.  
By using this topology $\mathcal O_{\rm sl}$, we give a new proof of Theorem A in the case where the sections are real projective spaces.  

\vspace{0.5truecm}

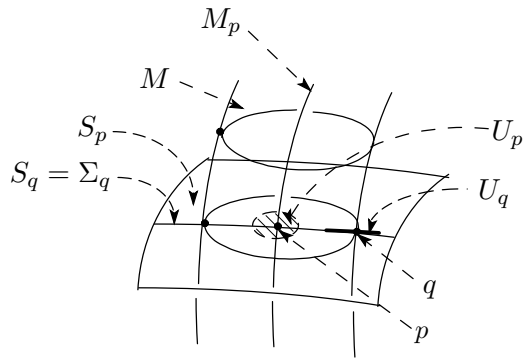
\centerline{
\unitlength 0.1in
\begin{picture}( 32.5300, 20.9000)( 13.8000,-30.3000)
%
\special{pn 8}%
\special{ar 3534 2492 756 864  3.2035570 4.1803242}%
%
\special{pn 8}%
\special{ar 4790 2578 756 864  3.2021026 4.1798978}%
%
\special{pn 8}%
\special{ar 2770 2898 2134 470  4.7187076 5.3540127}%
%
\special{pn 8}%
\special{ar 3132 2226 2134 468  4.7186868 5.3523394}%
%
\special{pn 20}%
\special{sh 1}%
\special{ar 3512 2106 10 10 0  6.28318530717959E+0000}%
\special{sh 1}%
\special{ar 3512 2106 10 10 0  6.28318530717959E+0000}%
%
\special{pn 8}%
\special{ar 4424 2554 526 1436  3.2015392 4.0680109}%
%
\special{pn 8}%
\special{ar 4424 2554 516 1426  2.9747157 3.1415927}%
%
\special{pn 8}%
\special{ar 3620 2492 524 1436  3.2015392 4.0691804}%
%
\special{pn 8}%
\special{ar 3602 2484 518 1426  2.9745992 3.1415927}%
%
\special{pn 8}%
\special{ar 4014 2510 526 1434  3.2014840 4.0688879}%
%
\special{pn 8}%
\special{ar 4000 2500 516 1426  2.9747157 3.1415927}%
%
\special{pn 8}%
\special{ar 2852 2444 2122 350  4.7187460 5.3534160}%
%
\special{pn 8}%
\special{pa 3358 1084}%
\special{pa 3682 1394}%
\special{da 0.070}%
\special{sh 1}%
\special{pa 3682 1394}%
\special{pa 3648 1334}%
\special{pa 3644 1358}%
\special{pa 3620 1362}%
\special{pa 3682 1394}%
\special{fp}%
%
\special{pn 8}%
\special{pa 3000 1350}%
\special{pa 3328 1488}%
\special{da 0.070}%
\special{sh 1}%
\special{pa 3328 1488}%
\special{pa 3274 1444}%
\special{pa 3280 1468}%
\special{pa 3260 1482}%
\special{pa 3328 1488}%
\special{fp}%
%
\special{pn 20}%
\special{sh 1}%
\special{ar 3926 2134 10 10 0  6.28318530717959E+0000}%
\special{sh 1}%
\special{ar 3926 2134 10 10 0  6.28318530717959E+0000}%
\put(26.4000,-19.2000){\makebox(0,0)[rb]{$S_q=\Sigma_q$}}%
\put(29.5000,-13.8000){\makebox(0,0)[rb]{$M$}}%
\put(33.2000,-11.1000){\makebox(0,0)[rb]{$M_p$}}%
\put(42.7000,-23.9000){\makebox(0,0)[lt]{$q$}}%
\put(26.2000,-16.8000){\makebox(0,0)[rb]{$S_p$}}%
%
\special{pn 8}%
\special{pa 3916 2144}%
\special{pa 3904 2174}%
\special{pa 3886 2200}%
\special{pa 3858 2218}%
\special{pa 3830 2232}%
\special{pa 3800 2244}%
\special{pa 3770 2254}%
\special{pa 3740 2262}%
\special{pa 3708 2266}%
\special{pa 3676 2272}%
\special{pa 3644 2274}%
\special{pa 3612 2276}%
\special{pa 3580 2278}%
\special{pa 3548 2276}%
\special{pa 3516 2276}%
\special{pa 3484 2272}%
\special{pa 3452 2270}%
\special{pa 3420 2266}%
\special{pa 3390 2262}%
\special{pa 3358 2254}%
\special{pa 3328 2246}%
\special{pa 3296 2236}%
\special{pa 3266 2226}%
\special{pa 3238 2212}%
\special{pa 3208 2198}%
\special{pa 3182 2182}%
\special{pa 3158 2160}%
\special{pa 3138 2134}%
\special{pa 3128 2104}%
\special{pa 3128 2096}%
\special{sp}%
%
\special{pn 8}%
\special{ar 3522 2088 392 142  3.1415927 4.6396779}%
%
\special{pn 8}%
\special{pa 3582 1954}%
\special{pa 3614 1956}%
\special{pa 3646 1960}%
\special{pa 3676 1966}%
\special{pa 3708 1970}%
\special{pa 3740 1978}%
\special{pa 3770 1984}%
\special{pa 3800 1994}%
\special{pa 3830 2006}%
\special{pa 3860 2018}%
\special{pa 3888 2034}%
\special{pa 3912 2054}%
\special{pa 3930 2080}%
\special{pa 3934 2112}%
\special{pa 3934 2116}%
\special{sp}%
%
\special{pn 8}%
\special{pa 4014 1668}%
\special{pa 4004 1698}%
\special{pa 3986 1724}%
\special{pa 3960 1742}%
\special{pa 3932 1760}%
\special{pa 3904 1772}%
\special{pa 3874 1784}%
\special{pa 3844 1792}%
\special{pa 3812 1798}%
\special{pa 3780 1804}%
\special{pa 3748 1808}%
\special{pa 3716 1810}%
\special{pa 3684 1810}%
\special{pa 3652 1810}%
\special{pa 3620 1810}%
\special{pa 3588 1808}%
\special{pa 3556 1806}%
\special{pa 3526 1800}%
\special{pa 3494 1794}%
\special{pa 3462 1790}%
\special{pa 3432 1780}%
\special{pa 3400 1772}%
\special{pa 3370 1762}%
\special{pa 3342 1748}%
\special{pa 3312 1734}%
\special{pa 3286 1718}%
\special{pa 3260 1698}%
\special{pa 3238 1676}%
\special{pa 3224 1646}%
\special{pa 3218 1616}%
\special{pa 3218 1616}%
\special{sp}%
%
\special{pn 8}%
\special{pa 3218 1640}%
\special{pa 3226 1610}%
\special{pa 3246 1584}%
\special{pa 3272 1566}%
\special{pa 3300 1550}%
\special{pa 3328 1536}%
\special{pa 3358 1526}%
\special{pa 3390 1518}%
\special{pa 3422 1512}%
\special{pa 3452 1508}%
\special{pa 3484 1502}%
\special{pa 3516 1500}%
\special{pa 3548 1498}%
\special{pa 3580 1498}%
\special{pa 3610 1496}%
\special{sp}%
%
\special{pn 8}%
\special{ar 3640 1650 374 144  4.8644791 6.2831853}%
\special{ar 3640 1650 374 144  0.0000000 0.0862056}%
\put(42.3000,-26.3000){\makebox(0,0)[lt]{{\small $p$}}}%
%
\special{pn 20}%
\special{ar 3110 2400 2208 286  5.0137088 5.1479380}%
%
\special{pn 8}%
\special{ar 3500 2100 120 70  0.0771195 0.7086984}%
\special{ar 3500 2100 120 70  1.0876458 1.7192247}%
\special{ar 3500 2100 120 70  2.0981721 2.7297511}%
\special{ar 3500 2100 120 70  3.1086984 3.7402774}%
\special{ar 3500 2100 120 70  4.1192247 4.7508037}%
\special{ar 3500 2100 120 70  5.1297511 5.7613300}%
\special{ar 3500 2100 120 70  6.1402774 6.3278923}%
%
\special{pn 4}%
\special{pa 3442 2162}%
\special{pa 3400 2120}%
\special{fp}%
\special{pa 3510 2170}%
\special{pa 3408 2068}%
\special{fp}%
\special{pa 3560 2160}%
\special{pa 3450 2050}%
\special{fp}%
\special{pa 3592 2132}%
\special{pa 3500 2040}%
\special{fp}%
\special{pa 3620 2100}%
\special{pa 3566 2046}%
\special{fp}%
%
\special{pn 8}%
\special{ar 2670 2060 300 210  4.7294522 4.9647463}%
\special{ar 2670 2060 300 210  5.1059228 5.3412169}%
\special{ar 2670 2060 300 210  5.4823933 5.7176875}%
\special{ar 2670 2060 300 210  5.8588639 6.0941581}%
\special{ar 2670 2060 300 210  6.2353345 6.2831853}%
%
\special{pn 8}%
\special{pa 2970 2070}%
\special{pa 2970 2090}%
\special{da 0.070}%
\special{sh 1}%
\special{pa 2970 2090}%
\special{pa 2990 2024}%
\special{pa 2970 2038}%
\special{pa 2950 2024}%
\special{pa 2970 2090}%
\special{fp}%
%
\special{pn 8}%
\special{ar 2650 1950 420 340  4.7123890 4.8702837}%
\special{ar 2650 1950 420 340  4.9650206 5.1229153}%
\special{ar 2650 1950 420 340  5.2176521 5.3755469}%
\special{ar 2650 1950 420 340  5.4702837 5.6281785}%
\special{ar 2650 1950 420 340  5.7229153 5.8808100}%
\special{ar 2650 1950 420 340  5.9755469 6.1334416}%
\special{ar 2650 1950 420 340  6.2281785 6.2441427}%
%
\special{pn 8}%
\special{pa 3070 1970}%
\special{pa 3070 1990}%
\special{da 0.070}%
\special{sh 1}%
\special{pa 3070 1990}%
\special{pa 3090 1924}%
\special{pa 3070 1938}%
\special{pa 3050 1924}%
\special{pa 3070 1990}%
\special{fp}%
%
\special{pn 8}%
\special{pa 4220 2380}%
\special{pa 3930 2150}%
\special{da 0.070}%
\special{sh 1}%
\special{pa 3930 2150}%
\special{pa 3970 2208}%
\special{pa 3972 2184}%
\special{pa 3996 2176}%
\special{pa 3930 2150}%
\special{fp}%
%
\special{pn 8}%
\special{pa 4000 2100}%
\special{pa 3990 2140}%
\special{da 0.070}%
\special{sh 1}%
\special{pa 3990 2140}%
\special{pa 4026 2080}%
\special{pa 4004 2088}%
\special{pa 3988 2070}%
\special{pa 3990 2140}%
\special{fp}%
%
\special{pn 8}%
\special{pa 3580 2050}%
\special{pa 3598 2024}%
\special{pa 3616 1996}%
\special{pa 3634 1972}%
\special{pa 3654 1946}%
\special{pa 3676 1922}%
\special{pa 3698 1900}%
\special{pa 3720 1878}%
\special{pa 3744 1854}%
\special{pa 3768 1834}%
\special{pa 3794 1814}%
\special{pa 3818 1794}%
\special{pa 3844 1776}%
\special{pa 3872 1758}%
\special{pa 3898 1742}%
\special{pa 3926 1726}%
\special{pa 3954 1710}%
\special{pa 3984 1698}%
\special{pa 4012 1684}%
\special{pa 4042 1670}%
\special{pa 4072 1660}%
\special{pa 4102 1648}%
\special{pa 4132 1640}%
\special{pa 4164 1630}%
\special{pa 4194 1622}%
\special{pa 4226 1616}%
\special{pa 4256 1608}%
\special{pa 4288 1604}%
\special{pa 4320 1600}%
\special{pa 4352 1596}%
\special{pa 4384 1592}%
\special{pa 4416 1590}%
\special{pa 4448 1588}%
\special{pa 4480 1588}%
\special{pa 4512 1588}%
\special{pa 4544 1590}%
\special{pa 4576 1592}%
\special{pa 4608 1594}%
\special{pa 4634 1598}%
\special{sp 0.070}%
%
\special{pn 8}%
\special{ar 4470 2140 490 230  3.4465467 3.6132134}%
\special{ar 4470 2140 490 230  3.7132134 3.8798800}%
\special{ar 4470 2140 490 230  3.9798800 4.1465467}%
\special{ar 4470 2140 490 230  4.2465467 4.4132134}%
\special{ar 4470 2140 490 230  4.5132134 4.6798800}%
\special{ar 4470 2140 490 230  4.7798800 4.7801583}%
\put(46.3000,-17.0000){\makebox(0,0)[lb]{$U_p$}}%
\put(45.6000,-20.1000){\makebox(0,0)[lb]{$U_q$}}%
\put(31.1000,-30.3000){\makebox(0,0)[lt]{$U_p,\,U_q\in\mathcal O_{\rm sl}$}}%
%
\special{pn 20}%
\special{sh 1}%
\special{ar 3130 2090 10 10 0  6.28318530717959E+0000}%
\special{sh 1}%
\special{ar 3130 2090 10 10 0  6.28318530717959E+0000}%
%
\special{pn 20}%
\special{sh 1}%
\special{ar 3210 1610 10 10 0  6.28318530717959E+0000}%
\special{sh 1}%
\special{ar 3210 1610 10 10 0  6.28318530717959E+0000}%
%
\special{pn 8}%
\special{pa 3570 2060}%
\special{pa 3580 2090}%
\special{da 0.070}%
\special{sh 1}%
\special{pa 3580 2090}%
\special{pa 3578 2020}%
\special{pa 3564 2040}%
\special{pa 3540 2034}%
\special{pa 3580 2090}%
\special{fp}%
%
\special{pn 8}%
\special{pa 4210 2640}%
\special{pa 3520 2130}%
\special{da 0.070}%
\special{sh 1}%
\special{pa 3520 2130}%
\special{pa 3562 2186}%
\special{pa 3564 2162}%
\special{pa 3586 2154}%
\special{pa 3520 2130}%
\special{fp}%
\end{picture}%
\hspace{2.2truecm}}

\vspace{0.5truecm}

\centerline{{\bf Figure  5$\,:\,$ Slice topology}}

\vspace{0.5truecm}

\noindent
{\it Proof of Theorem A (projective section-case)}\ \ 
Suppose that there exists an equifocal submanifold $M$ with real projective spatial section whose codimension is greater than two 
in an irreducible symmetric space $G/K$ of compact type and rank greater than one.  
Let $\pi:\widehat{G/K}\to G/K$ be the universal covering of the topological space $(G/K,\mathcal O_{\rm sl})$ 
and $\Gamma$ the deck transformation group of $\pi$.  Then we can determine uniquely the manifold structure and the Riemannian metric of $\widehat{G/K}$ 
such that $\pi$ is a Riemannian submersion onto the symmetric space $G/K$.  Let $\widehat{\mathfrak F}$ be the foliation consisting of the inverse images of 
the leaves of $\mathfrak F$ by $\pi$.  Then we can show that $\widehat{\mathfrak F}$ is a polar foliation with spherical section.  
As in the previous section, we can construct an irreducible topological Tits building of spherical type and rank greater than two for $\widehat{\mathfrak F}$ 
and hence $\widehat{G/K}$ is homothetic to $S^{n+r}(1)$.  It is clear that $\Gamma$ acts on $S^{n+r}(1)$ isometrically and freely.  Also, it is easy to show 
that $\Gamma$ is compact and connected as a subgroup of the isometry group of $\widehat{G/K}$.  
From these facts, we can derive that $\Gamma$ is isomorphic to the trivial group $\{e\}$, the circle group $U(1)$ or the special unitary group $SU(2)$.  
That is, $G/K$ is isometric to a sphere, a complex projective space or a quaternionic projective space.  This contradicts that $G/K$ is 
of rank greater than one.  Therefore, there does not exist the above equifocal submanifold.  \qed

\section{Proof of Theorem B} 
In this section, we prove Theorem B.  
For its purpose, we prepare some lemmas for rank one Lie triple systems.  
Let $G/K$ be an irreducible simply connected symmetric space of compact type and rank greater than one.  
Denote  by $r_{G/K}$ the rank of $G/K$.  
Let $\mathfrak g$ (resp. $\mathfrak k$) be the Lie algebra of $G$ (resp. $K$) and 
$\theta$ be the Cartan involution of $G$ with $({\rm Fix}\,\theta)_0\subset K\subset{\rm Fix}\,\theta$.  
Denote by the same symbol the involution of $\mathfrak g$ induced from $\theta$.  
Let $\mathfrak g=\mathfrak k\oplus\mathfrak p$ be the eigenspace decomposition of $\theta$.  
Note that $\mathfrak p$ is identified with the tangent space of $G/K$ at $eK$, where $e$ is the identity element of $G$.  
Denote by ${\rm Exp}$ the exponential map of $G/K$ at $eK$ and by $\exp$ the exponential map of $G$.  
Any totally geodesic submanifold in $G/K$ through $eK$ is given as the image of a Lie triple system of $\mathfrak p$ by ${\rm Exp}$ and 
it is a symmetric space.  Let $\mathfrak t$ be a Lie triple system of $\mathfrak p$.  If ${\rm Exp}(\mathfrak t)$ is a symmetric space of 
rank $r$, then we call $\mathfrak t$ (resp. ${\rm Exp}(\mathfrak t)$) a {\it rank $r$ Lie triple system} (resp. a {\it rank $r$ 
totally geodesic submanifold}).  Take a maximal abelian subspace $\mathfrak a$ of $\mathfrak p$.  Let $\mathfrak c$ be a Weyl domain in 
$\mathfrak a$.  Denote by $S(1)$ the unit sphere centered the origin in $\mathfrak p$.  Let ${\rm Ad}_G:G\to{\rm GL}(\mathfrak g)$ be 
the adjoint representation of $G$ and $\rho:K\to{\rm GL}(\mathfrak p)$ be the s-representation of $G/K$, that is, 
$\rho(k)(v):={\rm Ad}_G(k)(v)\,\,(k\in K,\,v\in\mathfrak p)$.  All $\rho(K)$-orbits meet the closure $\overline{\mathfrak c}$ of 
$\mathfrak c$ at the only one point, that is, $\overline{\mathfrak c}$ is the orbit space of the $\rho(K)$-action.  
The $\rho(K)$-orbits meeting $\mathfrak c$ are the principal orbits of the $\rho(K)$-action and they are full irreducible isoparametric 
submanifolds of codimension $r_{G/K}$ in $\mathfrak p$.  In particular, the $\rho(K)$-orbits meeting $\mathfrak c\cap S(1)$ are isoparametric 
submanifolds of codimension $r_{G/K}-1$ in $S(1)$.  Also, the $\rho(K)$-orbits meeting the boundary $\partial\mathfrak c$ of $\mathfrak c$ are 
the singular orbits of the $\rho(K)$-action and they are focal submanifolds of the principal orbits.  

\vspace{0.5truecm}

\noindent
{\bf Lemma 5.1.} {\sl Any rank one Lie triple system of $\mathfrak p$ is included by the cone 
$$C(\rho(K)(\vv_0)):=\{s\rho(k)(\vv_0)\,\vert\,k\in K,\,s\in{\mathbb R}\}$$
over the isotropy orbit $\rho(K)(\vv_0)$ for some $\vv_0\in\overline{\mathfrak c}\cap S(1)$.}

\vspace{0.5truecm}

\noindent
{\it Proof.} 
Set $\mathfrak k':=[\mathfrak t,\mathfrak t]$ and $\mathfrak g':=\mathfrak k'+\mathfrak t$.  
Since $\mathfrak t$ is the Lie triple system, we have 
$[\mathfrak k',\mathfrak k']\subset\mathfrak k',\,\,[\mathfrak g',\mathfrak g']\subset\mathfrak g'$ and 
$[\mathfrak k',\mathfrak t]\subset\mathfrak t$.  
Thus $\mathfrak g'$ and $\mathfrak k'$ are Lie subalgebras of $\mathfrak g$.  
Let $G'$ (resp. $K'$) be the connected Lie subgroup of $G$ with Lie algebra 
$\mathfrak g'$ (resp. $\mathfrak k'$) and ${\widetilde G}'$ (resp. ${\widetilde K}'$) the universal covering 
group of $G'$ (resp. $K'$).  Since $[\mathfrak k',\mathfrak t]\subset\mathfrak t$ and 
$[\mathfrak t,\mathfrak t]\subset\mathfrak k'$, ${\widetilde G}'/{\widetilde K}'$ is a simply connected 
symmetric space and it is the universal covering of ${\rm Exp}(\mathfrak t)$.  
Denote by $\rho'$ (resp. ${\widetilde{\rho}}'$) the isotropy representation of $G'/K'$ 
(resp. ${\widetilde G}'/{\widetilde K}'$).  
Take any unit vectors $v$ and $v'$ of $\mathfrak t$.  
Since ${\rm Exp}(\mathfrak t)$ is of rank one, both ${\rm Span}\{v\}$ and ${\rm Span}\{v'\}$ 
are a maximal abelian subspaces of $\mathfrak t$.  
Hence there exists $k'\in{\widetilde K}'$ with ${\widetilde{\rho}}'(k'){\rm Span}\{v\}={\rm Span}\{v'\}$.  
Since $\mathfrak g'\subset\mathfrak g$ and 
$\mathfrak k'=[\mathfrak t,\mathfrak t]\subset[\mathfrak p,\mathfrak p]\subset\mathfrak k$, 
$G'$ (resp. $K'$) is a Lie subgroup of $G$ (resp. $K$) and hence 
${\widetilde{\rho}}'({\widetilde K}')=\rho'(K')\subset\rho(K)$.  
Hence we obtain ${\widetilde{\rho}}'(k')\in\rho(K)$.  
The orbit $\rho(K)(v)(={\rm Ad}_G(K)(v))$ meets $\overline{\mathfrak c}$ 
at the only one point.  Let $\vv_0$ be this intersection point.  
The unit vectors $v$ and $v'$ belong to the orbit $\rho(K)(\vv_0)$.  
From the arbitrarinesses of $v$ and $v'$, it follows that $\mathfrak t$ is included by the cone 
$C(\rho(K)(\vv_0))$.  This completes the proof.  \qed

\vspace{0.5truecm}

We call $\vv_0\in\overline{\mathfrak c}$ in the statement of Lemma 5.1 the {\it characteristic direction} 
of $\mathfrak t$ (or ${\rm Exp}(\mathfrak t)$).  

\vspace{0.5truecm}

\noindent
{\bf Lemma 5.2.} {\sl Let $\Sigma$ and $\Sigma'$ be rank one totally geodesic submanifolds in $G/K$ through $eK$.  
If the characteristic directions of $\Sigma$ and $\Sigma'$ are distinct, then the component of $\Sigma\cap\Sigma'$ containing $eK$ is 
the one-point set $\{eK\}$.}

\vspace{0.5truecm}

\noindent
{\it Proof.} 
Let $\vv_0$ (resp. $v'_0$) be the characteristic direction of $\Sigma$ (resp. $\Sigma'$) and set 
$\mathfrak t:=T_{eK}\Sigma$ and $\mathfrak t':=T_{eK}\Sigma'$.  Since $\mathfrak t$ (resp. $\mathfrak t'$) is included by the cone 
$C(\rho(K)(\vv_0))$ (resp. $C(\rho(K)(\vv'_0))$) by Lemma 5.1 and $C(\rho(K)(\vv_0))\cap C(\rho(K)(\vv'_0))=\{0\}$ because of 
$\vv_0\not=\vv'_0$, we have $\mathfrak t\cap\mathfrak t'=\{0\}$.  This implies that the component of $\Sigma\cap\Sigma'$ containing $eK$ 
is the one-point set $\{eK\}$.  \qed

\vspace{0.5truecm}

Let $\triangle$ be the root system of $(G,K)$ with respect to $\mathfrak a$ and 
$\triangle_+$ be the positive root system of $\triangle$ under some lexicographic ordering of 
$\mathfrak a^{\ast}$.  Denote by $\mathfrak p_{\alpha}$ the root space for $\alpha\in\triangle_+$.  
Then we have the following root space decomposition:
$$\mathfrak p=\mathfrak a\oplus\left(\mathop{\oplus}_{\alpha\in\triangle_+}\mathfrak p_{\alpha}\right).$$
Fix $\vv_0\in\overline{\mathfrak c}\cap S(1)$.  
Denote by $A^{\vv_0}$ (resp. $h_{\vv_0}$) the shape tensor (resp. the second fundamental form) of 
the orbit $M^{\vv_0}:=\rho(K)(\vv_0)$ in $S(1)$.  
Set $\triangle_+^{\vv_0}:=\{\alpha\in\triangle_+\,\vert\,\alpha(\vv_0)=0\}$.  
Then the tangent space $T_{\vv_0}M^{\vv_0}$ of $M^{\vv_0}$ at $\vv_0$ is given by 
$$T_{\vv_0}M^{\vv_0}=\mathop{\oplus}_{\alpha\in\triangle_+\setminus\triangle_+^{\vv_0}}\mathfrak p_{\alpha}$$
and the normal space $T^{\perp}_{\vv_0}M^{\vv_0}$ of $M^{\vv_0}$ in $S(1)$ at $\vv_0$ is given by 
$$T^{\perp}_{\vv_0}M^{\vv_0}=(\mathfrak a\ominus{\rm Span}\{\vv_0\})\oplus\left(
\mathop{\oplus}_{\alpha\in\triangle_+^{\vv_0}}\mathfrak p_{\alpha}\right),$$
where $\mathfrak a\ominus{\rm Span}\{\vv_0\}$ means $\mathfrak a\cap({\rm Span}\{\vv_0\})^{\perp}$.  
Let $H_{\alpha}$ be the vector of $\mathfrak a$ with $\alpha(\cdot)=\langle H_{\alpha},\cdot\rangle$ and 
set $\displaystyle{e_{\alpha}:=\frac{H_{\alpha}}{\vert\vert H_{\alpha}\vert\vert}}$, where 
$\langle\,\,,\,\,\rangle$ is the inner product of $\mathfrak p$.  
The shape tensor $A^{\vv_0}$ satisfies 
$$A^{\vv_0}_{\xi}\vert_{\mathfrak p_{\alpha}}=-\frac{\alpha(\xi)}{\alpha(\vv_0)}{\rm id}\quad
(\alpha\in\triangle_+\setminus\triangle_+^{\vv_0},\,\,\,\xi\in T^{\perp}_{\vv_0}M^{\vv_0}).\leqno{(5.1)}$$

\vspace{0.5truecm}

By using Lemma 5.1 and $(5.1)$, we can derive the following fact.  

\vspace{0.5truecm}

\noindent
{\bf Proposition 5.3.} {\sl (i)\ If a vector subspace $\mathfrak s$ of $\mathfrak p$ is included by the cone $C(\rho(K)(\vv_0))$, then 
$${\rm dim}\,\mathfrak s\leq\mathop{\max}_{\alpha\in\triangle_+}(m_{\alpha}+m_{2\alpha})+1.$$

(ii)\ $\displaystyle{\frac{H_{\alpha}}{\|H_{\alpha}\|}}$ ($\alpha\in\triangle_+$) only can be chracteristic directions of rank one 
Lie triple systems.
}

\vspace{0.5truecm}

\noindent
{\it Proof.} Let $k_0$ be an element of $K$ with $\rho(k_0)(\vv_0)\in\mathfrak s\cap S(1)$.  
The set $\rho(k_0)^{-1}(\mathfrak s\cap S(1))$ is included by $M^{\vv_0}:=C(\rho(K)(\vv_0))$ by the assumption and 
it is totally geodesic in $S(1)$.  By using these facts and $(5.1)$, we obtain 
$$\begin{array}{l}
\displaystyle{T_{\vv_0}(\rho(k_0)^{-1}(\mathfrak s\cap S(1)))
\subset
\left\{\vw\in T_{\vv_0}M^{\vv_0}\,\vert\,A^{\vv_0}_{\xi}\vw=0\,\,(\forall\,\xi\in T^{\perp}_{\vv_0}M^{\vv_0})\right\}}\\
\displaystyle{=\left\{\left.\vw=\sum_{\alpha\in\triangle_+\setminus\triangle_+^{\vv_0}}\vw_{\alpha}\in\mathfrak p\,\right\vert\,
\sum_{\alpha\in\triangle_+\setminus\triangle_+^{\vv_0}}\frac{\alpha(\xi)}{\alpha(\vv_0)}\vw_{\alpha}\in{\rm Span}\{\vv_0\}\,\,
(\forall\,\xi\in T^{\perp}_{\vv_0}M^{\vv_0})\right\}}\\
\displaystyle{=\left\{\left.\vw=\sum_{\alpha\in\triangle_+\setminus\triangle_+^{\vv_0}}\vw_{\alpha}\in\mathfrak p\,\right\vert\,
\alpha(\xi)\vw_{\alpha}=0\,\,(\forall\,\xi\in\mathfrak a\ominus{\rm Span}\{\vv_0\},\,\,
\forall\,\alpha\in\triangle_+\setminus\triangle_+^{\vv_0})\right\}}\\
\displaystyle{\subset\bigoplus_{\alpha\in\triangle_+\setminus\triangle_+^{\vv_0}\,\,{\rm s.t.}\,\,H_{\alpha}\in{\rm Span}\{\vv_0\}}
\mathfrak p_{\alpha},}
\end{array}$$
($\vw_{\alpha}:$ the $\mathfrak p_{\alpha}$-component of $\vw$).  This implies that 
$$\rho(k_0)^{-1}(\mathfrak s)\subset{\rm Span}\{H_{\alpha}\}\oplus\mathfrak p_{\alpha}\oplus\mathfrak p_{2\alpha}\leqno{(5.2)}$$
holds for some $\alpha\in\triangle_+$.  
Therefore we can derive the statement (i).  

Let $\mathfrak t$ be a rank one Lie triple system and $\vv_0$ its characteristic direction.  
Since $\mathfrak t$ is included by the cone $C(\rho(K)(\vv_0))$ by Lemma 5.1, it follows from the above proof of the statement (i) that 
$$\rho(k_0)^{-1}(\mathfrak t)\subset{\rm Span}\{H_{\alpha}\}\oplus\mathfrak p_{\alpha}\oplus\mathfrak p_{2\alpha}\leqno{(5.3)}$$
holds for some $\alpha\in\triangle_+$.  This implies that the characteristic direction of $\mathfrak t$ equals to 
$\displaystyle{\frac{H_{\alpha}}{\|H_{\alpha}\|}}$.  Therefore, the statement (ii) is derived.  \qed

\vspace{0.5truecm}

Now we shall prove Theorem B by using these lemmas and Fact 1 stated in Introduction.  

\vspace{0.5truecm}

\noindent
{\it Proof of Theorem B.} Let $M$ be as in the statement of Theorem C.  Denote by $r$ the codimension of $M$.  
According to Fact 1 stated in Introduction, the section $\Sigma_x$ of $M$ through $x$ is a totally geodesic sphere (or real projective space) 
of constant curvature in $G/K$.  Take $p\in|(\mathcal S_M)_1|\setminus|\mathcal V_M|$ and the focal normal vector field $\widetilde{\vv}$ of $M$ 
such that $p:=\eta_{\widetilde{\vv}}(x)$.  
Denote by $\mathcal D_{\widetilde{\vv}}^{\perp}$ the orthogonal complementary distribution of the focal distribution 
$\mathcal D_{\widetilde{\vv}}$.  Take any $\vw\in(\mathcal D_{\widetilde{\vv}}^{\perp})_x$. Let $\alpha:(-\varepsilon,\varepsilon)\to M$ be 
a $C^{\infty}$-curve satisfying $\alpha(0)=x$ and $\alpha'(s)\in(\mathcal D_{\widetilde{\vv}}^{\perp})_{\alpha(s)}$ 
($s\in(-\varepsilon,\varepsilon)$).  Define a geodesic variation $\delta:[0,1]\times(-\varepsilon,\varepsilon)\to G/K$ by 
$$\delta(t,s):=\gamma_{\widetilde{\vv}_{\alpha(s)}}(t)\quad\,\,((t,s)\in[0,1]\times(-\varepsilon,\varepsilon)),$$
where $\gamma_{\widetilde{\vv}_{\alpha(s)}}$ is the normal geodesic of the direction $\widetilde{\vv}_{\alpha(s)}$.  
Set $\displaystyle{J:=\left.\frac{\partial\delta}{\partial s}\right|_{s=0}}$, which is the Jacobi field along the normal geodesic 
$\gamma_{\widetilde{\vv}_x}$ with $J(0)=\vw$ and $J'(0)=-A_{\widetilde{\vv}_x}\vw$.  
Hence we have 
$$J(s)=P_{\gamma_{\widetilde{\vv}_x}|_{[0,s]}}\left(\cos\left(s\sqrt{\widetilde R(\widetilde{\vv}_x)}\right)(\vw)
-\frac{\sin\left(s\sqrt{\widetilde R(\widetilde{\vv}_x)}\right)}{\sqrt{\widetilde R(\widetilde{\vv}_x)}}(A_{\widetilde{\vv}_x}\vw)\right),$$
where $P_{\gamma_{\widetilde{\vv}_x}|_{[0,s]}}$ is the parallel translation along $\gamma_{\widetilde{\vv}_x}|_{[0,s]}$, 
$\widetilde R(\widetilde{\vv}_x)$ is the normal Jacobi operator $\widetilde R(\cdot,\widetilde{\vv}_x)\widetilde{\vv}_x$ for $\widetilde{\vv}_x$ 
($\widetilde R:$ the cuvature tensor of $G/K$) and $\displaystyle{\cos\left(s\sqrt{\widetilde R(\widetilde{\vv}_x)}\right)}$ 
(resp. $\displaystyle{\frac{\sin\left(s\sqrt{\widetilde R(\widetilde{\vv}_x)}\right)}{\sqrt{\widetilde R(\widetilde{\vv}_x)}}}$) is 
defined by 
$$\begin{array}{l}
\displaystyle{\cos\left(s\sqrt{\widetilde R(\widetilde{\vv}_x)}\right):=\sum_{k=0}^{\infty}\frac{(-1)^k}{(2k)!}\,
\widetilde R(\widetilde{\vv}_x)^k}\\
\displaystyle{\left({\rm resp.}\,\,\,\,\frac{\sin\left(s\sqrt{\widetilde R(\widetilde{\vv}_x)}\right)}{\sqrt{\widetilde R(\widetilde{\vv}_x)}}
:=\sum_{k=0}^{\infty}\frac{(-1)^k}{(2k+1)!}\,\widetilde R(\widetilde{\vv}_x)^k\,\,\right).}
\end{array}
\leqno{(5.3)}$$
Here we used a general description of Jacobi fields in a symmetric space (see Section 3 of [TT] or (1.2) of [K1]).  
Since $\Sigma_x$ is totally geodesic in $G/K$, $\widetilde R(\vv_x)$ preserves $T_x\Sigma_x$ invariantly.  
Hence it preserves $T_xM$ invariantly.  From this fact, we can derive $J(1)\in P_{\gamma_{\widetilde{\vv}_x}|_{[0,1]}}(T_xM)$.  
Also, we have $J(1)=(d\eta_{\widetilde{\vv}})_x(\vw)$.  Therefore, from the arbitrariness of $\vw$, we obtain 
$$T_pF_{\widetilde{\vv}}=(d\eta_{\widetilde{\vv}})_x(T_xM)\subset P_{\gamma_{\widetilde{\vv}_x}|_{[0,1]}}(T_xM).$$
On the other hand, since $\Sigma_x$ is totally geodesic in $G/K$, we have 
$$P_{\gamma_{\widetilde{\vv}_x}|_{[0,1]}}(T_x^{\perp}M)=P_{\gamma_{\widetilde{\vv}_x}|_{[0,1]}}(T_x\Sigma_x)=T_p\Sigma_x.$$
From these facts, it follows that $\Sigma_x$ meets $F_{\widetilde{\vv}}$ orthogonally at $p$.  
Similarly, we can show that, for any $y\in L^{\widetilde{\vv}}_x$, $\Sigma_y$ meets $F_{\widetilde{\vv}}$ orthogonally at $p$.  
Hence we can derive $\sum\limits_{y\in L^{\widetilde{\vv}}_x}T_p\Sigma_y\subset T_p^{\perp}F_{\widetilde{\vv}}$.  
Furthermore, since $L^{\widetilde{\vv}}_x$ is a fibre of the focal map $\eta_{\widetilde{\vv}}:M\to F_{\widetilde{\vv}}$, we can show that 
$\sum\limits_{y\in L^{\widetilde{\vv}}_x}T_p\Sigma_y=T_p^{\perp}F_{\widetilde{\vv}}$.  
This implies that the normal umbrella $S_p:=\exp^{\perp}(T_p^{\perp}F_{\widetilde{\vv}})$ is equal to 
$\displaystyle{\mathop{\cup}_{y\in L^{\widetilde{\vv}}_x}\Sigma_y}$.  
Without loss of generality, we may assume that $p=eK$.  
Since $\Sigma_y$'s $(y\in L^{\widetilde{\vv}}_x$) are totally geodesic rank one symmetric spaces, $T_p\Sigma_y$'s $(y\in L^{\widetilde{\vv}}_x$) 
are rank one Lie triple systems of $\mathfrak p=T_{eK}(G/K)$.  Also, since $p\in|(\mathcal S_M)_1|\setminus|\mathcal V_M|$, 
we can show $\displaystyle{{\rm dim}\left(\mathop{\cap}_{y\in L^{\widetilde{\vv}}_x}\Sigma_y\right)=1}$.  
According to Lemma 5.2, it follows from these facts that the Lie triple systems $T_p\Sigma_y$'s $(y\in L^{\widetilde{\vv}}_x$) have the same 
characteristic direction.  Denote by $\vv_0$ this characteristic direction.  According to Lemma 5.1, $T_p\Sigma_y$'s 
$(y\in L^{\widetilde{\vv}}_x$) are included by the cone $C(\rho(K)(\vv_0))$ and hence $T_p^{\perp}F_{\widetilde{\vv}}$ is also included by 
the cone $C(\rho(K)(\vv_0))$.  Therefore, by (i) of Proposition 5.3, we have 
$${\rm dim}\,T_p^{\perp}F_{\widetilde{\vv}}\leq\mathop{\max}_{\alpha\in\triangle_+}(m_{\alpha}+m_{2\alpha})+1.$$
On the other hand, since $p\in|(\mathcal S_M)_1|\setminus|\mathcal V_M|$ again, we have ${\rm dim}\,L^{\widetilde{\vv}}_x\geq r-1$.  
From these facts, we can derive 
$$r\leq\left(\mathop{\max}_{\alpha\in\triangle_+}(m_{\alpha}+m_{2\alpha})+1\right)-(r-1),$$
that is, 
$$r\leq\left[\frac{1}{2}\,\mathop{\max}_{\alpha\in\triangle_+}(m_{\alpha}+m_{2\alpha})\right]+1.$$
\qed

\vspace{0.5truecm}

{\small 
\rightline{Department of Mathematics, Faculty of Science}
\rightline{Tokyo University of Science, 1-3 Kagurazaka}
\rightline{Shinjuku-ku, Tokyo 162-8601 Japan}
\rightline{(koike@rs.tus.ac.jp)}
}

\end{document}